\documentclass[twoside,11pt,reqno]{amsart}
\usepackage{amsmath,amssymb,amscd,mathrsfs}

\makeatletter

\@addtoreset{equation}{section}

\def\thesubsection{\S\thesection-\alph{subsection}}
\def\Point#1{\addtocounter{subsection}{1}\vspace{2mm}
\noindent\thesubsection. {\bf #1.}\def\@currentlabel{\thesubsection}}

\newtheorem{Proposition}{Proposition}[section]
\newtheorem{Lemma}[Proposition]{Lemma}
\newtheorem{Theorem}[Proposition]{Theorem}
\newtheorem{Corollary}[Proposition]{Corollary}

\newbox\squ  % box character for ends of proofs
\setbox\squ=\hbox{\vrule width.6pt
   \vbox{\hrule height.6pt width.4em\kern1ex
         \hrule height.6pt}%
   \vrule width.6pt}
\def\Col{\operatorname{Col}}

\def\gr{\operatorname{gr}}
\def\op{{\operatorname{op}}}
\def\F{\operatorname{F}}

\def\C{{\mathbb C}}

\def\Z{{\mathbb Z}}

\def\0{{\bar 0}}
\def\1{{\bar 1}}

\def\hom{{\operatorname{Hom}}}

\def\End{{\operatorname{End}}}

\def\MP{{\mathcal{M}}}
\def\OP{{\mathcal{P}}}

\def\bc{\operatorname{\text{\boldmath$c$}}}
\def\bi{\operatorname{\text{\boldmath$i$}}}
\def\bj{\operatorname{\text{\boldmath$j$}}}
\def\bq{{\operatorname{\text{\boldmath$q$}}}}
\def\bt{{\operatorname{\text{\boldmath$t$}}}}
\def\blambda{{\operatorname{\text{\boldmath$\lambda$}}}}
\def\bnu{{\operatorname{\text{\boldmath$\nu$}}}}

\def\eps{{\varepsilon}}
\def\phi{{\varphi}}

\newdimen\hoogte    \hoogte=12pt    
\newdimen\breedte   \breedte=14pt  
\newdimen\dikte     \dikte=0.5pt 
\newenvironment{Young}{\begingroup
       \def\vr{\vrule height0.89\hoogte width\dikte depth 0.2\hoogte}
       \def\fbox##1{\vbox{\offinterlineskip
                    \hrule height\dikte
                    \hbox to \breedte{\vr\hfill##1\hfill\vr}
                    \hrule height\dikte}}
       \vbox\bgroup \offinterlineskip \tabskip=-\dikte \lineskip=-\dikte
            \halign\bgroup &\fbox{##\unskip}\unskip  \crcr }
       {\egroup\egroup\endgroup}
\def\diagram#1{\relax\ifmmode\vcenter{\,\begin{Young}#1\end{Young}\,}\else%
              $\vcenter{\,\begin{Young}#1\end{Young}\,}$\fi}

\begin{document}

\begin{abstract}
We prove that the center of each degenerate cyclotomic Hecke algebra
associated to the complex reflection group of type $B_d(l)$
consists of symmetric polynomials in 
its commuting generators.
The classification of the blocks of the degenerate cyclotomic
Hecke algebras is an easy consequence.
We then deduce that the center of an integral block of parabolic
category $\mathcal O$ for the Lie algebra $\mathfrak{gl}_n(\C)$
is generated by the center of its universal enveloping algebra.
\end{abstract}

\title[Degenerate cyclotomic Hecke algebras]{\sc Centers of Degenerate Cyclotomic Hecke Algebras and  Parabolic Category $\mathcal O$}
\author{\sc Jonathan Brundan}

%\address{Department of Mathematics, University of Oregon, Eugene, OR 97403, USA}
%\email{brundan@uoregon.edu}

\thanks{2000 {\it Mathematics Subject Classification}: 20C08, 17B20.} 
\thanks{Research supported in part by NSF grant no. DMS-0139019.}

\maketitle

\section{Introduction}

Let $R$ be a fixed commutative ground ring.
Recall from \cite{D} that the {\em degenerate affine Hecke algebra} $H_d$
is the $R$-algebra which is equal as an $R$-module to the 
tensor product $R[x_1,\dots,x_d] \otimes_R R S_d$ 
of the polynomial algebra
$R[x_1,\dots,x_d]$ and the group algebra $R S_d$ of the symmetric group $S_d$.
Multiplication is defined so that $R[x_1,\dots,x_d]$ 
(identified with the subspace $R[x_1,\dots,x_d] \otimes 1$)
and $R S_d$ (identified with the subspace $1 \otimes R S_d$) are subalgebras, 
and in addition
\begin{align*}
s_i x_{i+1} &= x_{i}s_i + 1,\\
s_i x_j &= x_j s_i \qquad\qquad\qquad(j \neq i,i+1),
\end{align*}
where $s_i$ denotes the basic transposition $(i\,i\!+\!1) \in S_d$.
It is known 
by 
\iffalse
the degenerate analogue 
of a result of Bernstein
\fi\cite[Theorem 6.5]{L}
that the center $Z(H_d)$ of $H_d$ consists of all symmetric polynomials
in the (algebraically independent) generators $x_1,\dots,x_d$.

Given
in addition a monic polynomial
$f(x) = x^l + c_1 x^{l-1} +\cdots + c_l \in R[x]$
of degree $l \geq 1$, 
the {\em degenerate cyclotomic Hecke algebra} 
$H_d^f$ is the quotient of $H_d$ by the two-sided ideal generated
by $f(x_1)$. We refer to $l$ here as the {\em level}.
Since we seldom
mention $H_d$ itself again, it should not cause confusion
to also use the notation $x_1,\dots,x_d$ for the canonical images of
the polynomial generators of $H_d$ in the quotient $H_d^f$. 
For example, if $f(x) = x$ then $H_d^f$ can be identified simply
with the group algebra $RS_d$, and under this identification
we have that
$$
x_i = \sum_{j=1}^{i-1} (j\,i) \:\in R S_d,
$$
the {\em Jucys-Murphy elements}.
In this case, it has long been known (see \cite{J} or \cite[1.9]{M})
that the center of $R S_d$
again consists of all symmetric polynomials in 
$x_1,\dots,x_d$, though of course these generators are no longer algebraically
independent.
In other words, the canonical homomorphism $H_d \twoheadrightarrow R S_d$
maps $Z(H_d)$ surjectively onto
$Z(R S_d)$.
Our first result proves the analogous statement for the
quotient map $H_d \twoheadrightarrow H_d^f$ in general.

\vspace{2mm}
\noindent
{\bf Theorem 1.} {\em The center of $H_d^f$
consists of all symmetric polynomials
in $x_1,\dots,x_d$.
Moreover, $Z(H_d^f)$ is free as an $R$-module
with an explicit basis parametrized by all $l$-multipartitions of $d$.}
\vspace{2mm}

For the first application, specialize to the case that $R = F$ is an algebraically closed field.
We say that two irreducible modules $L$ and $L'$ 
{\em belong to the same block} if they are {linked}
by a chain
$L = L_0,L_1,\dots,L_n = L'$ 
of irreducible modules 
such that there is a non-split extension between
$L_{i-1}$ and $L_i$ for each $i=1,\dots,n$.
For modules over a finite dimensional algebra like $H_d^f$, 
this is equivalent to 
the property that $L$ and $L'$ have the same central character.
So, on combining Theorem 1 with the existing theory,
we obtain the classification of 
the blocks of the degenerate cyclotomic Hecke algebras.\footnote{
In an earlier version of this article,
we also explained how to deduce 
the classification of blocks of the degenerate
affine Hecke algebra $H_d$ from Theorem 1. However, Iain Gordon has pointed
out that this follows immediately by a general result of M\"uller
\cite[III.9.2]{BG},
since $H_d$ is finite as a module over its center.
}
The conclusion is exactly as claimed in Grojnowski's unpublished note
\cite{Gnote}; see 
$\S$4 below for the precise statement.
Unfortunately, as has been pointed out by Anton Cox,
the argument given there 
is incomplete, so this corrects an error in the literature.
\iffalse
The same mistake was made in the twisted case 
in \cite[Corollary 8.13]{BKcrystal}, so at the end of this article we also discuss
briefly how to extend Theorem 1 to this situation, replacing the degenerate
cyclotomic
Hecke algebra with the degenerate cyclotomic Hecke-Clifford superalgebra.
\fi
Actually, \cite{Gnote} was mainly concerned with
cyclotomic Hecke algebras (not their rational degenerations). 
For these, it has also long been expected that the center consists of all 
symmetric polynomials in the Jucys-Murphy elements, but we still do not know how
to prove this. Nevertheless, Lyle and Mathas \cite{LM} have recently 
managed to solve the problem of classifying the blocks of the 
cyclotomic Hecke algebras too, by a quite different method.

Now we further specialize to the case that
$F = \C$.
Let $\mu = (\mu_1,\dots,\mu_l)$ be an $l$-tuple of positive integers summing to $n$.
Let $\mathfrak{g} = \mathfrak{gl}_n(\C)$
and let $\mathfrak{p}$ be the standard parabolic subalgebra with 
block diagonal Levi
subalgebra $\mathfrak{h} = \mathfrak{gl}_{\mu_1}(\C) \oplus\cdots\oplus
\mathfrak{gl}_{\mu_l}(\C)$. 
Let $\mathcal O^\mu$ 
be the category of all {\em finitely generated} $\mathfrak{g}$-modules
which are {\em locally finite} as $\mathfrak{p}$-modules and 
{\em integrable} as $\mathfrak{h}$-modules, 
i.e. they lift 
to rational representations of $H = GL_{\mu_1}(\C) \times\cdots\times GL_{\mu_l}(\C)$. This is the usual parabolic analogue of the BGG category
$\mathcal O$, except that we are only allowing modules with integral 
weights/central characters.
The category $\mathcal O^\mu$ decomposes as
$$
\mathcal O^\mu = \bigoplus_\nu \mathcal O^\mu_\nu
$$
where the direct sum is over integral central
characters $\nu:Z(\mathfrak{g}) \rightarrow \C$ 
of the universal enveloping algebra $U(\mathfrak{g})$,
and $\mathcal O^\mu_\nu$ is the full subcategory of
$\mathcal O^\mu$ consisting of modules with generalized
central character $\nu$.

The next result, also ultimately a consequence of Theorem 1, is an essential
ingredient in \cite{Bpre,St}, which give
quite different (and independent) proofs of a conjecture
of Khovanov \cite[Conjecture 3]{Kh}.
Recall that the center $Z(\mathcal C)$
of an additive category $\mathcal C$ is the commutative ring consisting of all
natural transformations from the identity functor to itself.
For example, if $\mathcal C$ is the category of finite dimensional 
modules over a finite dimensional
algebra $C$, then $Z(\mathcal C)$ is canonically isomorphic to the center
of the algebra $C$ itself.

\vspace{2mm}
\noindent
{\bf Theorem 2.}
{\em
For any integral central character $\nu$,
the natural map
$$
m^\mu_\nu: 
Z(\mathfrak{g}) \rightarrow Z(\mathcal O^\mu_\nu)
$$
sending $z \in Z(\mathfrak{g})$ to the natural transformation
defined by left multiplication by $z$
is a surjective algebra homomorphism.
Moreover, the dimension of $Z(\mathcal O^\mu_\nu)$ is the same as
the number of isomorphism classes of irreducible modules
in $\mathcal O^\mu_\nu$.
}
\vspace{2mm}

The category $\mathcal O^\mu_\nu$ is equivalent to the category
of finite dimensional 
modules over a finite dimensional algebra,
e.g. one can take endomorphism algebra
of a minimal projective generator.
Hence, two irreducible modules in $\mathcal O^\mu_\nu$
belong to the same block if and only if 
they have the same central character with respect to 
$Z(\mathcal O^\mu_\nu)$.
By definition, all irreducible modules in $\mathcal O^\mu_\nu$
have the same central character with respect to $Z(\mathfrak{g})$.
So Theorem 2 implies that all the irreducible modules in
$\mathcal O^\mu_\nu$ belong to the same block.
This proves that the above decomposition of 
$\mathcal O^\mu$ (defined by central 
characters) coincides with its decomposition into blocks 
in the usual sense
(defined by linkage classes of irreducible modules).
For regular central characters in arbitrary type, this is
a known consequence of some results of Deodhar combined with 
the Kazhan-Lusztig conjecture, 
but for singular central characters even in type $A$
this was an open problem. According to Boe, the same coincidence
is expected in types $D$ and $E$, but there are 
counterexamples in non-simply-laced types.

\iffalse
In \cite[$\S$6]{Kh}, Khovanov has formulated 
several interesting conjectures about the
structure of the centers of the categories $\mathcal O(\mu,\theta)$,
some of which have already been proved in \cite[$\S$5.2]{MS}.
The most intriguing of these asserts that, 
for any regular central character 
$\theta$, the center of $\mathcal O(\mu,\theta)$ is canonically 
isomorphic to the (complexified) cohomology algebra of the Springer fiber 
associated to a nilpotent matrix of Jordan type $\mu$.
In a forthcoming article \cite{Bpre} I hope to use Theorem 2 
to prove this conjecture. 
Moreover, we will give an explicit
presentation for $Z(\mathcal O(\mu,\theta))$ 
for singular central characters too.
By Schur-Weyl duality for higher levels, this 
also yields presentations for the centers of arbitrary blocks of the
degenerate cyclotomic Hecke algebras.
\fi

The remainder of the article is organized as follows.
There is a natural filtration on the algebra $H_d^f$
with respect to which
the associated graded algebra $\gr H_d^f$ is 
the twisted tensor
product of the level $l$
truncated polynomial algebra $R[x_1,\dots,x_d] / (x_1^l,\dots,x_d^l)$ by the group algebra $R S_d$ of the symmetric group.
In section 2, 
we compute the center of this associated graded algebra directly,
giving the crucial upper bound on the size of $Z(H_d^f)$
since we obviously have that $\gr Z(H_d^f) \subseteq Z(\gr H_d^f)$.
There are then several different 
ways to show that this upper bound is actually attained.
The approach followed in section 3 is to simply
write down enough linearly independent 
central elements in $H_d^f$.
This has the advantage of yielding at the same time an explicit
basis for $Z(H_d^f)$ which is a generalization of Murphy's basis 
for $Z(RS_n)$ constructed in the proof of \cite[1.9]{M}.
In section 4, we discuss the classification of the blocks of $H_d^f$
in more detail. In particular we compute the dimension of the center
of each block, refining Theorem 1 which gives the dimension of the
center of the whole algebra.
Finally in section 5 we deduce the results about parabolic category
$\mathcal O$ by exploiting the 
Schur-Weyl 
duality for higher levels from \cite{schur}, which reduces many questions
about the category $\mathcal O^\mu$ to the 
degenerate cyclotomic Hecke algebras $H_d^f$ for
$f(x) = (x-\mu_1) \cdots (x-\mu_l)$.

\vspace{2mm}
\noindent
{\bf Acknowledgements.} Thanks
 to Brian Boe, Alexander Kleshchev and Victor Ostrik
for helpful discussions, and Iain Gordon for pointing out M\"uller's theorem.

\section{The center of the associated graded algebra}

We fix an integer $l \geq 1$
and a commutative ring $R$.
Let $R_l[x_1,\dots,x_d]$ denote the level $l$ truncated polynomial algebra,
that is, the quotient of the polynomial algebra $R[x_1,\dots,x_d]$ by the
relations $x_1^l = \cdots = x_d^l = 0$. The symmetric group $S_d$
acts on $R_l[x_1,\dots,x_d]$ by algebra automorphisms so that
$w \cdot x_i = x_{wi}$ for each $i$ and $w \in S_d$.
We view the resulting twisted tensor product algebra
$R_l[x_1,\dots,x_d] \,{\scriptstyle{\rtimes\!\!\!\!\!\bigcirc}}\, R S_d$
as a graded algebra with each $x_i$ in degree $1$ and all elements
of $S_d$ in degree $0$.
The goal in this section is to compute the center of this algebra explicitly.
We remark that the algebra 
$R_l[x_1,\dots,x_d] \,{\scriptstyle{\rtimes\!\!\!\!\!\bigcirc}}\, R S_d$
can be viewed as a degeneration of the group algebra
$R(C_l \wr S_d)$ of the wreath product of the symmetric group and 
the cyclic group of order $l$.
It is well known that the conjugacy classes of $C_l \wr S_d$
are parametrized by certain multipartitions; see \cite[p.170]{Mac}
or \cite{wang}. 
With this in mind the results in  this section
should not be too surprising.

Let $Q_d$ denote the centralizer of $R_l[x_1,\dots,x_d]$ in
$R_l[x_1,\dots,x_d] \,{\scriptstyle{\rtimes\!\!\!\!\!\bigcirc}}\, R S_d$.
The symmetric group $S_d$ acts 
on $Q_d$ by conjugation, i.e.
$w \cdot z = wzw^{-1}$.
It is obvious that the center 
of $R_l[x_1,\dots,x_d] \,{\scriptstyle{\rtimes\!\!\!\!\!\bigcirc}}\, R S_d$ is just the set of fixed points:
$$
Z(R_l[x_1,\dots,x_d] \,{\scriptstyle{\rtimes\!\!\!\!\!\bigcirc}}\, R S_d) = Q_d^{S_d}.
$$
We are going first to describe an explicit basis for $Q_d$
from which it will be easy to determine the $S_d$-fixed points,
hence the center.

For $r \geq 0$ and any 
set $I = \{i_1,\dots,i_a\}$ of $a$ distinct numbers chosen from
$\{1,\dots,d\}$, let
$$
h_r(I) = h_r(i_1,\dots,i_a)
:= 
\!\!\!\!\sum_{\substack{0 \leq r_1,\dots,r_a < l \\ r_1+\cdots+r_a=(a-1)(l-1)+r}} x_{i_1}^{r_1} \cdots x_{i_a}^{r_a} \:\in R_l[x_1,\dots,x_d],
$$
the $((a-1)(l-1)+r)$th complete symmetric function in the variables
$x_{i_1},\dots,x_{i_a}$.
By the pigeonhole principle, $h_r(I)$ is zero if $r \geq l$,
and moreover $h_{l-1}(I) = x_{i_1}^{l-1} \cdots x_{i_a}^{l-1}$.

\begin{Lemma}\label{hl}
Let $I, J$ be any two subsets of $\{1,\dots,d\}$ 
with $c = |I \cap J| > 0$.
For any $r,s \geq 0$, we have that
$h_r(I) h_s(J) = l^{c-1} h_{r+s+(c-1)(l-1)}(I \cup J)$.
\end{Lemma}

\begin{proof}
Suppose first that $I = \{i_1,\dots,i_a,k\}$ and
$J = \{j_1,\dots,j_b,k\}$ with $I \cap J = \{k\}$.
Then we have that
\begin{align*}
h_r(I) h_s(J) &= \sum_{\substack{r_1,\dots,r_a \\ s_1,\dots,s_b}}
x_{i_1}^{r_1} \cdots x_{i_a}^{r_a} x_k^{a(l-1)+r-r_1-\cdots-r_a}
x_{j_1}^{s_1} \cdots x_{j_b}^{s_b} x_k^{b(l-1)+s-s_1-\cdots-s_b}\\
&= \sum_{\substack{r_1,\dots,r_a \\ s_1,\dots,s_b}}
x_{i_1}^{r_1} \cdots x_{i_a}^{r_a}
x_{j_1}^{s_1} \cdots x_{j_b}^{s_b}
 x_k^{(a+b)(l-1)+r+s-r_1-\cdots-r_a-s_1-\cdots-s_b}\\
&= 
h_{r+s}(I \cup J).
\end{align*}
This proves the lemma in the case $c=1$.
Next we take $i \neq j$ and note that
\begin{align*}
h_0(i,j) h_0(i,j)
&=
\sum_{r,s} x_i^{r} x_j^{l-1-r} 
x_j^{s} x_i^{l-1-s}\\
&=
\sum_{r,s} x_i^{l-1+r-s} x_j^{l-1+s-r}
= l x_i^{l-1} x_j^{l-1}=
l h_{l-1}(i,j).
\end{align*}
Finally take $I = \{i_1,\dots,i_a,k_1,\dots,k_c\}$,
$J = \{j_1,\dots,j_b,k_1,\dots,k_c\}$ for $c \geq 2$
and assume that $I \cap J = \{k_1,\dots,k_c\}$.
Using the preceeding two formulae, we get that
\begin{align*}
h_r(I) h_s(J)
&= 
h_r(i_1,\dots,i_a,k_1)
h_0(k_1,k_2)
\cdots 
h_0(k_{c-1},k_c) \\
&\qquad\qquad\qquad\,\quad\times
h_0(k_1,k_2)
\cdots 
h_0(k_{c-1},k_c)
h_s(j_1,\dots,j_b,k_1)\\
&=
l^{c-1}
h_{r+s+(c-1)(l-1)}(i_1,\dots,i_a,j_1,\dots,j_b,k_1,\dots,k_c)
\\
&= l^{c-1} h_{r+s+(c-1)(l-1)}(I \cup J).
\end{align*}
This is what we wanted.
\end{proof}

Now let $A = (i_1 \, \cdots \, i_a)$ be an $a$-cycle in $S_d$.
Write $h_r(A)$ for $h_r(i_1,\dots,i_a)$.
Given another cycle $B = (j_1\,\cdots\,j_b)$,
write
$A \cup B$ and $A \cap B$ for the {\em sets} 
$\{i_1,\dots,i_a\} \cup \{j_1,\dots,j_b\}$ and
$\{i_1,\dots,i_a\} \cap \{j_1,\dots,j_b\}$, respectively.
Let 
$$
A^{(r)} := h_r(A) A\:\in
R_l[x_1,\dots,x_d] \,{\scriptstyle{\rtimes\!\!\!\!\!\bigcirc}}\, R S_d,
$$
which we call a cycle of {\em color $r$}.
As before, we have that $A^{(r)} = 0$ for $r \geq l$,
so we need only consider colors from the set $\{0,1,\dots,l-1\}$.
In the case of $1$-cycles, we have that
$(i)^{(r)} = x_i^r$, so $1$-cycles of color $0$ are trivial.

\begin{Lemma}\label{cand}
Let $r, s \geq 0$ be colors and $A$ and $B$ be cycles in $S_d$.
Let $c = |A \cap B|$.
\begin{itemize}
\item[(i)]
If $c = 0$ then $A^{(r)} B^{(s)} = B^{(s)} A^{(r)}$, i.e.
disjoint colored cycles commute.
\item[(ii)]
If $c = 1$ (in which case the product $AB$ is a single cycle)
then
$$
A^{(r)} B^{(s)} = (AB)^{(r+s)}.
$$
\item[(iii)] If $c \geq 2$ then
$$
A^{(r)} B^{(s)} = \delta_{r+s,0} l^{c-1} h_{(c-1)(l-1)}(A \cup B) AB.
$$
\end{itemize}
\end{Lemma}

\begin{proof}
(i) Obvious.

(ii) Say $A = (i_1\,\cdots\, i_a\,k)$ and $B = (j_1\,\cdots\,j_b\,k)$.
By Lemma~\ref{hl}, we have that
\begin{align*}
A^{(r)} B^{(s)} &= h_r(i_1,\dots,i_a,k) A
h_s( j_1,\dots, j_b, k) B\\
&=
h_r(i_1,\dots,i_a,k) h_s(j_1,\dots,j_b,i_1)
AB\\
&=
h_{r+s}(i_1,\dots,i_a,j_1,\dots,j_b,k) AB
= (AB)^{(r+s)}
\end{align*}
as required.

(iii) Arguing exactly as in (ii), we get that
$$
A^{(r)} B^{(s)} = l^{c-1} h_{r+s+(c-1)(l-1)}(A \cup B) AB.
$$
Now observe that $h_{r+s+(c-1)(l-1)}(A \cup B)$ is zero unless $r+s=0$.
\end{proof}

Now we are going to consider products of colored cycles.
Using Lemma~\ref{cand}, it is easy to see that any such product is
either zero or else it can be rewritten 
as some power of $l$ times a {\em product of disjoint colored cycles},
meaning a product $A_1^{(r_1)} \cdots A_m^{(r_m)}$
where $A_1,\dots,A_m$ are disjoint cycles in $S_d$ 
and $0 \leq r_1,\dots,r_m < l$ are some colors.
Moreover, 
two such products of disjoint colored cycles are equal if and only if
one can be obtained from the other by reordering the 
disjoint colored cycles and adding/removing
some $1$-cycles of color $0$.
For example,
\begin{align*}
(1\,2\,3)^{(4)} (7\,9\,2)^{(1)} &= ((1\,2\,3)(7\,9\,2))^{(5)}
=
(1\,2\,7\,9\,3)^{(5)}
=
(1\,2\,7\,9\,3)^{(5)} (4)^{(0)},\\
(1\,2\,3)^{(4)} (7\,9\,2\,1)^{(1)} &= 0 = (1\,2\,3)^{(0)} (7\,9\,3\,2\,1)^{(0)}
\text{ (assuming $l > 1$)},\\
(1\,2\,3)^{(0)} (7\,9\,2\,1)^{(0)} &= 
l(x_1x_2x_3x_7x_9)^{l-1} (1\,2\,3)(7\,9\,2\,1)
=
l(1\,7\,9\,3)^{(l-1)} (2)^{(l-1)}.
\end{align*}

\begin{Theorem}\label{basd}
The set of all products of disjoint colored cycles
is a basis for $Q_d$. In particular, it is a free $R$-module
of rank 
$$
\sum
\frac{d!}{r_1! r_2!\cdots} \left(\frac{l}{1}\right)^{r_1} \left(\frac{l}{2}\right)^{r_2} \cdots
$$ summing over all partitions 
$(1^{r_1} 2^{r_2}\cdots)$ of $d$.
\end{Theorem}

\begin{proof}
Observe by applying Lemma~\ref{cand} with $B^{(s)} = (i)^{(1)} =  x_i$
that every colored cycle $A^{(r)}$ belongs to the algebra $Q_d$.
Hence all products of colored cycles belong to $Q_d$.
Moreover, 
any product of colored cycles is a linear combination of products of 
disjoint colored
cycles, and the
set of all products of disjoint colored cycles is linearly independent.
It just remains to show that $Q_d$ is spanned
by all products of colored cycles.

Suppose to start with that $A = (i_1\,\cdots\,i_a)$ is an $a$-cycle 
in $S_d$ and that
$$
z=\sum_{\substack{0 \leq r_1,\dots,r_a < l \\ r_1+\cdots+r_a=k}} c_{r_1,\dots,r_a} x_{i_1}^{r_1}\cdots x_{i_a}^{r_a} A
$$
is a non-zero homogeneous element of $Q_d$
of degree $k \geq 0$ for some coefficients $c_{r_1,\dots,r_a} \in R$.
We claim that $z$ is a scalar multiple of
$A^{(r)}$ for some $0 \leq r < l$ (in which case $k=(a-1)(l-1)+r$).
To see this,
equating coefficients of $x_{i_1}^{r_1} \cdots x_{i_j}^{r_j+1} \cdots
x_{i_a}^{r_a} A$
in the equation $x_{i_j} z = z x_{i_j}$ gives that
$$
c_{r_1,\dots,r_j,r_{j+1},\dots,r_a} 
= 
c_{r_1,\dots,r_j+1,r_{j+1}-1,\dots,r_a}
$$
whenever $r_j < l-1$ for some $j=1,\dots,a-1$, interpreting the right hand side as zero in case $r_{j+1}=0$.
If $k = (a-1)(l-1)+r$ for $0 \leq r < l$, we deduce from this that
all the coefficents $c_{r_1,\dots,r_a}$ are equal
to 
$c_{l-1,\dots,l-1,r}$, hence $z$ is a scalar multiple of $A^{(r)}$.
Otherwise, we can write $k = m(l-1)+r$ for some
$m \leq a-2$ and $0 \leq r < l-1$ and get that all the coefficients
$c_{r_1,\dots,r_a}$ are equal to
$c_{l-1,\dots,l-1,r,0,0,\dots,0} = c_{l-1,\dots,l-1,r+1,-1,0,\dots,0} = 0$,
contradicting the assumption that $z \neq 0$.

Now take an element $f w \in Q_d$ for $w \in S_d$ and a homogeneous polynomial
$f \in R_l[x_1,\dots,x_d]$.
Write $w = A_1 \cdots A_m$ as a product of disjoint cycles, none of which are
$1$-cycles.
We show by induction on $m$ that $fw$ is a linear combination of products
of colored cycles. The base case $m=0$ is clear as then $w = 1$.
For the induction step, suppose that $m \geq 1$ and
$A_m = (i_1\,\cdots\,i_a)$. Let $I = \{i_1,\dots,i_a\}$ and
$J = \{1,\dots,d\}\setminus I$.
We can write $f = \sum_{s=1}^t f_s g_s$ for
homogeneous polynomials $f_1,\dots,f_t \in R_l[x_i\mid i \in I]$
and linearly independent homogeneous polynomials
$g_1,\dots,g_t
\in R_l[x_j\mid j \in J]$.
Equating coefficients of $g_sA_1 \cdots A_{m-1}$ 
in the equations $x_i fw = fw x_i$
for each $i \in I$, we deduce that each $f_s A_m$ belongs to $Q_d$.
Hence by the previous paragraph each $f_s A_m$ is a scalar multiple
of $A_m^{(r)}$ for some $0 \leq r < l-1$.
This shows that $fw = \sum_{r=0}^{l-1} h_r A_1\cdots A_{m-1} A_m^{(r)}$
for homogeneous polynomials $h_r \in R_l[x_j\mid j \in J]$.
Equating coefficients of $A_m^{(r)}$ in the equations $x_j fw = fwx_j$
for each $j \in J$, we deduce that each
$h_r A_1 \cdots A_{m-1}$ belongs to $Q_d$. Hence by the induction hypothesis
each $h_r A_1 \cdots A_{m-1}$ is a linear combination of products of
colored cycles. Hence $fw$ is too.

Finally take an arbitrary homogeneous element 
$\sum_{w \in S_d} f_w w \in Q_d$, for polynomials
$f_w \in R_l[x_1\dots,x_d]$.
We have for each $i$ that $\sum_{w\in S_d} x_i f_w w = 
\sum_{w \in S_d} x_{wi} f_w w$.
Equating coefficients 
gives that $x_i f_w w = x_{wi} f_w w = f_w w x_i$
for each $i$ and $w$. Hence each $f_w w$ belongs to $Q_d$.
So by the previous paragraph each $f_w w$ is a linear combination of
products of colored cycles. This completes the proof.
\end{proof}

For a partition $\lambda = (\lambda_1 \geq \lambda_2 \geq \cdots)$ 
we write $|\lambda|$ for $\lambda_1+\lambda_2+\cdots$ and
$\ell(\lambda)$ for its {\em length}, that is, the number of non-zero parts.
By an {\em $l$-multipartition of $d$} 
we mean a tuple $\blambda = (\lambda^{(1)},\lambda^{(2)},\dots,\lambda^{(l)})$
of partitions such that 
$|\lambda^{(1)}|+\cdots+|\lambda^{(l)}| = d$.
Let $\MP_d(l)$ denote the set of all
$l$-multipartitions of $d$.
Given a product $z = A_1^{(r_1)} \cdots A_m^{(r_m)}$ 
of disjoint colored cycles in $Q_d$, where each $A_i$ is an $a_i$-cycle,
we can add extra $1$-cycles of color $0$ if necessary to assume
that $a_1+\cdots+a_m = d$.
Define the 
{\em cycle type} of $z$ to be the {$l$-multipartition}
$\blambda = (\lambda^{(1)},\dots,\lambda^{(l)})$ of $d$ defined by
declaring that $\lambda^{(r)}$ is the partition whose parts consist of
all the $a_i$ such that $r_i = r-1$. 
For $\blambda = (\lambda^{(1)},\dots,\lambda^{(l)})
\in \MP_d(l)$, let $z_d(\blambda)$
denote the sum of all 
products of disjoint colored cycles in $Q_d$ of cycle type
$\blambda$.

\begin{Theorem}\label{ct}
The elements $\{z_d(\blambda)\mid \blambda \in \MP_d(l)\}$ form a basis for
the center of 
$R_l[x_1,\dots,x_d] \,{\scriptstyle{\rtimes\!\!\!\!\!\bigcirc}}\, R S_d$.
In particular, 
$Z(R_l[x_1,\dots,x_d] \,{\scriptstyle{\rtimes\!\!\!\!\!\bigcirc}}\, R S_d)$ 
is a free $R$-module of rank 
$|\MP_d(l)|$.
\end{Theorem}

\begin{proof}
As we remarked at the beginning of the section, the center of
$R_l[x_1,\dots,x_d] \,{\scriptstyle{\rtimes\!\!\!\!\!\bigcirc}}\, R S_d$
is the set of fixed points of $S_d$ on $Q_d$.
Given a colored cycle $A^{(r)} = (i_1\,\cdots\,i_a)$ and $w \in S_d$,
we have that
$$
w \cdot A^{(r)} = (w \cdot A)^{(r)} = (w i_1\,\cdots\,w i_a)^{(r)}.
$$
So the action of $S_d$ on $Q_d$ is the linear action induced by 
a permutation action on
the basis from Theorem~\ref{basd}. It just remains to observe that
two products of disjoint colored cycles lie in the same 
$S_d$-orbit if and only if they have the same cycle type,
and the $z_d(\blambda)$'s are simply the orbit sums.
\end{proof}

\begin{Corollary}\label{cc}
If $d!$ is invertible in $R$ then
the center of 
$R_l[x_1,\dots,x_d] \,{\scriptstyle{\rtimes\!\!\!\!\!\bigcirc}}\, R S_d$
is generated by the elements
$$
z_d(a^{(r)}) 
:= \sum_{\text{all $a$-cycles }A\in S_d} A^{(r)}
$$
for all $0 \leq r < l$ and $1 \leq a \leq d$.
\end{Corollary}

\begin{proof}
Take a 
multipartition $\blambda = (\lambda^{(1)},\dots,\lambda^{(l)}) \in \MP_d(l)$.
Consider the product of the elements $z_d(a^{(r-1)})$ over
all $r = 1,\dots,l$ and all non-zero parts $a$ of $\lambda^{(r)}$.
It gives an invertible scalar multiple of $z_d(\blambda)$ modulo 
lower terms.
\end{proof}

In the remainder of the section, we are going to 
construct another basis for 
$Z(R_l[x_1,\dots,x_d] \,{\scriptstyle{\rtimes\!\!\!\!\!\bigcirc}}\, R S_d)$
which is a generalization of the basis 
for the center of 
$RS_n$ constructed by Murphy in the proof of \cite[1.9]{M}.
Given $k \geq 0$ and $1 \leq i \leq d$,
write $k = (a-1)l+r$ for $a \geq 1$ and $0 \leq r < l$, then set
$$
y_i(k) := 
\sum_{\substack{1 \leq i_1,\dots,i_{a-1} < i \\ i_1,\dots,i_{a-1}\text{ distinct}}}
(i_1\,\cdots\,i_{a-1}\,i)^{(r)},
$$
an element of degree $(a-1)(l-1)+r$.
For example, 
$y_i(k) = 0$ if $k \geq il$ and
$y_i(r) = (i)^{(r)} = x_i^r$ for $0 \leq r < l$.
Particularly important, we have that
$$
y_i(l) = \sum_{j=1}^{i-1} (j\,i)^{(0)},
$$
which we call the $i$th {\em colored Jucys-Murphy element}.

\begin{Lemma}\label{cjm}
For any $1 \leq i \leq d$ and $p \geq 0$, we have that
$$
y_i(l)^p = y_i(pl) + (*)
$$ where $(*)$ is a linear combination of products
$A_1^{(l-1)} \cdots A_m^{(l-1)}$
for disjoint cycles $A_1,\dots,A_m$ in $S_i$ such that $A_1$ involves $i$
and $|A_1 \cup \cdots \cup A_m| \leq p$.
\end{Lemma}

\begin{proof}
Induction exercise 
using Lemma~\ref{cand}.
\end{proof}

\iffalse
\begin{Corollary}
For any $1 \leq i \leq d$, $p \geq 0$ and $0 < r < l$, 
we have that $y_i(l)^p y_i(r) = y_i(pl+r)$.
\end{Corollary}
\fi

For a partition $\lambda = (\lambda_1 \geq \lambda_2 \geq \cdots)$,
let $\lambda/ l := (\lfloor \lambda_1 / l \rfloor \geq \lfloor \lambda_2 / l \rfloor \geq \cdots)$.
We are going to use partitions belonging to the
set
$$
\OP_d(l) = \{\lambda \mid \ell(\lambda) + |\lambda / l | \leq d\}
$$
to parametrize our new basis.
Note to start with 
that $|\OP_d(l)| = |\MP_d(l)|$, so this set is of the right size. Indeed, there is a bijection
$$
\varphi:\MP_d(l) \rightarrow \OP_d(l),
$$
defined as follows.
Suppose that 
$\blambda = (\lambda^{(1)},\dots,\lambda^{(l)}) \in \MP_d(l)$ 
where 
$\lambda^{(r)} = (\lambda^{(r)}_1 \geq \cdots \geq \lambda^{(r)}_{m_r} > 0)$.
Then $\phi(\blambda)$ denotes
the ordinary partition with parts
$(\lambda^{(r)}_i-1)l + r-1$ for all $1 \leq r \leq l$ and $1 \leq i \leq m_r$.
It is easy to see that $\phi(\blambda)$ belongs to $\OP_d(l)$.
Conversely, given $\mu = (\mu_1 \geq \mu_2 \geq \cdots) \in \OP_d(l)$,
there is a unique multipartition $\blambda = (\lambda^{(1)}, \dots,\lambda^{(l)})$ of $d$
such that $\varphi(\blambda) = \mu$:
the parts of $\lambda^{(r)}$
are the numbers  $\lfloor\mu_i / l \rfloor+1$ for all $i=1,\dots,d-|\mu/l|$ such that 
$\mu_i \equiv r-1 \pmod{l}$.
Hence $\phi$ is indeed a bijection.

Since every element $\mu$ of $\OP_d(l)$ is of length at most $d$,
it can be thought
of simply as a $d$-tuple of integers.
Given two $d$-tuples $\mu = (\mu_1,\dots,\mu_d)$ and
$\nu = (\nu_1,\dots,\nu_d)$ we write $\mu \sim \nu$
if one is obtained from the other by permuting the entries.
For $\mu = (\mu_1,\dots,\mu_d) \in \OP_d(l)$, define
$$
m_d(\mu) := \sum_{\nu \sim \mu} y_1(\nu_1) \cdots y_d(\nu_d).
$$
This is a homogeneous element of $Q_d$ of degree
$|\mu| - |\mu/l|$. 

\begin{Theorem}\label{bacon}
The elements $\{m_d(\mu)\mid \mu \in \OP_d(l)\}$ form a
basis for the center of $R_l[x_1,\dots,x_d] \,{\scriptstyle{\rtimes\!\!\!\!\!\bigcirc}}\, R S_d$.
\end{Theorem}

\begin{proof}
Let us first check that $m_d(\mu)$ belongs to 
$Z(R_l[x_1,\dots,x_d] \,{\scriptstyle{\rtimes\!\!\!\!\!\bigcirc}}\, R S_d)$.
We just need to check it commutes with each basic
transposition $(i\,i\!+\!1)$.
Obviously, $(i\,i\!+\!1)$ commutes with $y_j(m)$ 
if $j \neq i,i+1$. Therefore it suffices to show for each $i=1,\dots,d-1$
and $k, m \geq 0$ that $(i\,i\!+\!1)$ commutes with both the elements
$y_i(k) y_{i+1}(k)$
and
$y_i(k) y_{i+1}(m) + y_i(m) y_{i+1}(k)$.
For the first case, write $k = (a-1)l+r$ as usual.
We have that
$$
y_i(k) y_{i+1}(k) = 
\sum_{\substack{1 \leq i_1,\dots,i_{a-1} < i \\ 1 \leq j_1,\dots,j_{a-1} < i+1}}
 (i_1\,\cdots\, i_{a-1}\, i)^{(r)}
 (j_1\,\cdots\, j_{a-1}\, i\!+\!1)^{(r)}
$$
where the sum is over distinct $i_1,\dots,i_{a-1}$
and distinct $j_1,\dots,j_{a-1}$.
We split this sum into two pieces:
$$
\sum_{\substack{1 \leq i_1,\dots,i_{a-1} < i \\ 1 \leq j_1,\dots,j_{a-1} < i}}
 (i_1\,\cdots\, i_{a-1}\, i)^{(r)}
 (j_1\,\cdots\, j_{a-1}\, i\!+\!1)^{(r)}
$$
which clearly commutes with $(i\,i\!+\!1)$, and
$$
\sum_{b=1}^{a-1}
\sum_{\substack{1 \leq i_1,\dots,i_{a-1} < i \\ 1 \leq j_1,\dots,j_{a-1} < i+1 \\ j_b = i}}
 (i_1\,\cdots\, i_{a-1}\, i)^{(r)}
 (j_1\,\cdots\, j_{b-1} \, i\, j_{b+1}\,\cdots\,j_{a-1}\, i\!+\!1)^{(r)}
$$
which also commutes with $(i\,i\!+\!1)$ by an application of Lemma~\ref{cand}.
The second case is similar.

Now we compare the $m_d(\mu)$'s with the basis from 
Theorem~\ref{ct}. For any $\blambda \in \MP_d(l)$
define $\#\blambda$ to be $(d-z)$ where $z$ is the number of
parts of $\lambda^{(1)}$ that equal $1$.
We claim for $\blambda \in \MP_d(l)$ with $\phi(\blambda) = \mu$
that
$m_d(\mu) = z_d(\blambda) + (*)$
where $(*)$ is a linear combination of $z_d(\bnu)$'s
for $\bnu \in \MP_d(l)$ with $\#\bnu < \#\blambda$.
The theorem clearly follows from this claim and Theorem~\ref{ct}.
To prove the claim,
let $\mu = (\mu_1 \geq \cdots \geq \mu_h > 0)$ and write each
$\mu_i$ as $(a_i-1)l+r_i$ as usual,
so $\#\blambda = a_1+\cdots+a_h$.
By definition, $m_d(\mu)$ is a sum of products of colored cycles
of the form $z=A_1^{(r_1)} \cdots A_h^{(r_h)}$ 
where each $A_i$ is an $a_i$-cycle.
If $A_1,\dots,A_h$ happen to be disjoint cycles then $z$
is of cycle type $\blambda$.
Otherwise, $|A_1\cup\cdots\cup A_h| < \#\blambda$
so using Lemma~\ref{cand} we can rewrite $z$ as a linear combination of
products of disjoint colored cycles of cycle type
$\bnu \in \MP_d(l)$ with $\#\bnu < \#\blambda$.
Combined with the first paragraph and Theorem~\ref{ct}, 
this shows that $m_d(\mu) = c z_d(\blambda) + (*)$ for some $c$.
Finally, to show that $c=1$, consider the coefficient of 
one particular product of disjoint colored cycles of cycle type
$\blambda$ in the expansion of $m_d(\mu)$.
\end{proof}

\section{The center of $H_d^f$}

We are ready to tackle the problem of computing the center
of the degenerate cyclotomic Hecke algebra
$H_d^f$, where  $f(x) = x^l + c_1 x^{l-1} + \cdots + c_l \in R[x]$ 
is a monic polynomial of degree $l$.
Define a filtration 
$$
\F_0 H_d^f \subseteq \F_1 H_d^f \subseteq \F_2 H_d^f \subseteq \cdots
$$
of the algebra $H_d^f$ by declaring that $\F_r H_d^f$ is spanned
by all $x_{i_1} \cdots x_{i_s} w$ for  $0 \leq s \leq r, 
1 \leq i_1,\dots,i_s \leq d$ and $w \in S_d$.
So each $x_i$ is in filtered degree $1$ and each $w \in S_d$
is in filtered degree $0$.
Given an element $z \in \F_r H_d^f$, we write
$\gr_r z$ for its canonical image in the $r$th graded component
$\gr_r H_d^f = \F_r H_d^f / \F_{r-1} H_d^f$ of the associated graded
algebra $\gr H_d^f = \bigoplus_{r \geq 0} \gr_r H_d^f$.
By the PBW theorem for degenerate cyclotomic Hecke algebras
\cite[Lemma 3.5]{schur},
this associated graded algebra $\gr H_d^f$ can be identified with
the twisted tensor product
$R_l[x_1,\dots,x_d] \,{\scriptstyle{\rtimes\!\!\!\!\!\bigcirc}}\, R S_d$
so that $\gr_1 x_i$ is identified with
$x_i \in 
R_l[x_1,\dots,x_d] \,{\scriptstyle{\rtimes\!\!\!\!\!\bigcirc}}\, R S_d$
and $\gr_0 w$ is identified with 
$w \in 
R_l[x_1,\dots,x_d] \,{\scriptstyle{\rtimes\!\!\!\!\!\bigcirc}}\, R S_d$.
To avoid confusion, we reserve the notations
$x_i^r$ and $s_i$ from now on for the elements of $H_d^f$,
always using the alternate notations
$(i)^{(r)}$ and $(i\,i\!+\!1)$ for the corresponding elements of
the associated graded algebra 
$R_l[x_1,\dots,x_d] \,{\scriptstyle{\rtimes\!\!\!\!\!\bigcirc}}\, R S_d$.

Given an $R$-submodule $V$ of $H_d^f$,
we can consider the {\em induced filtration} on $V$
defined by setting $\F_r V := V \cap \F_r H_d^f$.
The associated graded module $\gr V$ is canonically identified with
an $R$-submodule of 
$R_l[x_1,\dots,x_d] \,{\scriptstyle{\rtimes\!\!\!\!\!\bigcirc}}\, R S_d$,
and for two submodules we have that $V = V'$ if and only if $\gr V = \gr V'$.
Note also that 
$$
\gr Z(H_d^f) \subseteq 
Z(R_l[x_1,\dots,x_d] \,{\scriptstyle{\rtimes\!\!\!\!\!\bigcirc}}\, R S_d).
$$
Hence if we can find elements $z_1 \in \F_{i_1} 
Z(H_d^f),\dots,z_m \in \F_{i_m} Z(H_d^f)$
with the property that $\gr_{i_1} z_1,\dots, \gr_{i_m} z_m$ 
is a basis for
$Z(R_l[x_1,\dots,x_d] \,{\scriptstyle{\rtimes\!\!\!\!\!\bigcirc}}\, R S_d)$,
then it follows immediately that $z_1,\dots,z_m$ also is a basis
for $Z(H_d^f)$. This is exactly what we are going to do.
Recall the elements $y_i(k)$ and $m_d(\mu)$ of 
$R_l[x_1,\dots,x_d] \,{\scriptstyle{\rtimes\!\!\!\!\!\bigcirc}}\, R S_d$
from the previous section.

\begin{Lemma}\label{wwagto}
Assume that $1 \leq i \leq d$ and $k = (a-1)l+r$ for some
$a \geq 1$ and $0 \leq r < l$.
Then we have that $x_i^k \in \F_{(a-1)(l-1)+r} H_d^f$ and
$$
\gr_{(a-1)(l-1)+r} x_i^k
=
\left\{ 
\begin{array}{ll}
y_i(k)+(*)&\text{if $r = 0$,}\\
y_i(k)&\text{if $r > 0$,}
\end{array}\right.
$$
where
$(*)$ denotes a linear combination of 
products of disjoint colored
cycles of the form $A_1^{(l-1)} \cdots A_m^{(l-1)}$ 
such that $i \in A_1 \cup \cdots \cup A_m \subseteq \{1,\dots,i\}$
and $|A_1 \cup \cdots \cup A_m| \leq a-1$.
\end{Lemma}

\begin{proof}
Assume to start with that $k=l$, i.e. $a=2,r=0$.
We prove the lemma in this case by induction on $i=1,\dots,d$.
For the base case,
we have that $x_1^l = -c_1 x_1^{l-1}-\cdots - c_l$,
so it is in filtered degree $(l-1)$ and $\gr_{l-1} x_1^l
= -c_1 (1)^{(l-1)} = y_1(l)-c_1(1)^{(l-1)}$.
For the induction step, we have by the relations that
$$
x_{i+1}^l = s_i x_i^l s_i + \sum_{t=0}^{l-1} x_i^t x_{i+1}^{l-1-t} s_i.
$$
Hence by induction we get that
$x_{i+1}^l$ is in filtered degree $(l-1)$ and
\begin{align*}
\gr_{l-1} x_{i+1}^l &= (i\,i\!+\!1) (y_i(l) - c_1 (i)^{(l-1)}) (i\,i\!+\!1)
+ (i\,i\!+\!1)^{(0)}\\
 &= y_{i+1}(l)  -c_1 (i+1)^{(l-1)}
\end{align*}
as we wanted.

Now assume that $k=(a-1)l$ for any $a \geq 1$, i.e. the case when $r=0$.
By the previous paragraph, we have that 
$x_i^k = (x_i^l)^{a-1}$ is in filtered degree $(a-1)(l-1)$
and
$$
\gr_{(a-1)(l-1)} x_i^k = (y_i(l) -c_1 (i)^{(l-1)})^{a-1}.
$$
By Lemma~\ref{cjm} this equals $y_i(k) + (*)$ where $(*)$
is a linear combination of 
products of disjoint colored
cycles of the form $A_1^{(l-1)} \cdots A_m^{(l-1)}$ 
such that $i \in A_1 \cup \cdots \cup A_m \subseteq \{1,\dots,i\}$
and $|A_1 \cup \cdots \cup A_m| \leq a-1$.

Finally assume that $k=(a-1)l+r$ for $0 < r < l$.
Writing $x_i^k = (x_i^{(a-1)l}) (x_i^r)$ and using the previous paragraph
and Lemma~\ref{cand} gives the desired conclusion in this case.
\end{proof}

For any $d$-tuple $\mu = (\mu_1,\dots,\mu_d)$ of non-negative integers,
let
$$
p_d(\mu) := \sum_{\nu \sim \mu} x_1^{\nu_1} \cdots x_d^{\nu_d}\:\in H_d^f.
$$
Since this is a symmetric polynomial in $x_1,\dots,x_d$, it is
automatically central.
Theorem 1 from the introduction is a consequence of the following 
more precise result.

\begin{Theorem}\label{pr}
For $\mu \in \OP_d(l)$, we have that
$p_d(\mu) \in \F_{r} Z(H_d^f)$ 
where $r = |\mu|-|\mu/l|$.
Moreover, 
$\gr_r p_d(\mu) = m_d(\mu) + (*)$
where $(*)$ is a linear combination of $m_d(\nu)$'s 
for $\nu \in \OP_d(l)$ with $|\nu/l|+\ell(\nu) < |\mu/l|+\ell(\mu)$.
Hence, $\gr Z(H_d^f) = Z(R_l[x_1,\dots,x_d] \,{\scriptstyle{\rtimes\!\!\!\!\!\bigcirc}}\, R S_d)$ and 
the elements
$$
\left\{p_d(\mu)\mid\mu \in \OP_d(l)\right\}
$$
form a basis for $Z(H_d^f)$. In particular, $Z(H_d^f)$ is a free $R$-module
of rank equal to the number of $l$-multipartitions of $d$.
\end{Theorem}

\begin{proof}
Recall the bijection $\phi:\MP_d(l) \rightarrow \OP_d(l)$ and also the notation
$\#\blambda$ from the last paragraph of the proof of Theorem~\ref{bacon}.
We showed there for $\blambda \in \MP_d(l)$ with $\phi(\blambda) = \mu$ that
$m_d(\mu) = z_d(\blambda) + (*)$
where $(*)$ is a linear combination of $z_d(\bnu)$'s with $\#\bnu < \#\blambda$.
Note $\#\blambda = |\mu / l| + \ell(\mu)$.
So we get from this also that
$z_d(\blambda) = m_d(\mu) + (*)$
where $(*)$ is a linear combination of $m_d(\nu)$'s for $\nu \in \OP_d(l)$
with $|\nu/l|+\ell(\nu) < |\mu/l|+\ell(\mu)$.

Now, by Lemma~\ref{wwagto} and the definitions, $p_d(\mu)$ is in filtered degree
$r = |\mu|-|\mu/l|$ and moreover
$\gr_r p_d(\mu) = m_d(\mu) + (*)$
where $(*)$ is a linear combination of products of disjoint colored cycles 
$A_1^{(r_1)}\cdots A_m^{(r_m)}$ such that 
$|A_1\cup\cdots\cup A_m| < |\mu/l|+\ell(\mu)$.
Since $\gr_r p_d(\mu)$ is central, it follows by Lemma~\ref{cand} and Theorem~\ref{ct}
that $(*)$ can be rewritten as a linear combination of $z_d(\bnu)$'s with
$\#\bnu < \#\blambda$. Hence by the first paragraph it is also a linear combination
of $m_d(\nu)$'s with $|\nu/l|+\ell(\nu) < |\mu/l|+\ell(\mu)$.
This proves that the elements
$$
\left\{\gr_{|\mu|-|\mu/l|} p_d(\mu)\mid \mu \in \OP_d(l)\right\}
$$
form a basis for $Z(R_l[x_1,\dots,x_d] \,{\scriptstyle{\rtimes\!\!\!\!\!\bigcirc}}\, R S_d)$.
Now the theorem follows by the general principles discussed just
before Lemma~\ref{wwagto}.
\end{proof}

\begin{Corollary}\label{ps}
If $d!$ is invertible in $R$, then
the center of $H_d^f$ is generated by
the power sums $x_1^r+\cdots+x_d^r$
for $1 \leq r \leq d$.
\end{Corollary}

\begin{proof}
Under the assumption on $R$,
it is well known that every symmetric polynomial in variables
$x_1,\dots,x_d$ can be expressed as a polynomial
in the first $r$ power sums.
\end{proof}

\section{The blocks of $H_d^f$}

In this section, we replace the ground ring $R$ with a
ground field $F$ such that we can factor
$f(x) = (x-q_1)\cdots (x-{q_l})$ for $q_1,\dots,q_l \in F$.
We will denote the algebra $H_d^f$ instead by $H_d^{\bq}$
where $\bq = (q_1,\dots,q_l) \in F^l$.
We point out that $F$ is a splitting field for the algebra $H_d^{\bq}$;
one way to see this is to check that
the construction of the irreducible 
$H_d^{\bq}$-modules over the algebraic closure of $F$
from \cite[$\S$5.4]{Kbook} already makes sense over $F$ itself.
Theorem~1 just proved shows in particular that
the dimension of 
$Z(H_d^\bq)$ is equal to the number of $l$-multipartitions of $d$.
The goal in this section is to refine this statement by 
computing the dimensions of the centers
of the individual blocks.

Before we can even formulate the result, 
we need an explicit combinatorial parametrization of 
the blocks, or equivalently, the central characters of $H_d^{\bq}$.
This is a well known consequence of Theorem 1.
To start with we recall the classification of central characters
of $H_d$ itself following \cite[$\S$4.2]{Kbook}.
Given a tuple $\bi = (i_1,\dots,i_d) \in F^d$, write
$$
\chi({\bi}): Z(H_d) \rightarrow F
$$
for the central character mapping a symmetric polynomial
$f(x_1,\dots,x_d)$ to $f(i_1,\dots,i_d)$.
Clearly, 
$\chi({\bi}) = \chi({\bj})$ if and only if $\bi \sim \bj$,
so this gives a parametrization of central characters of $H_d$
by the set $X_d$ of $\sim$-equivalence classes in $F^d$.
Now we pass to the quotient $H_d^{\bq}$ of $H_d$.
Since $Z(H_d)$ maps surjectively onto $Z(H_d^{\bq})$ by Theorem 1,
the set of all central characters of $H_d^{\bq}$
is naturally parametrized by the subset
$$
X_d^{\bq} = \{\bi \in X_d\mid\chi({\bi}):Z(H_d) \rightarrow F
\text{ factors through the quotient }Z(H_d^{\bq})\}
$$
of $X_d$.
To complete the classification
of blocks of $H_d^{\bq}$,
it just remains to describe this subset $X_d^{\bq}$ combinatorially.

To do this, we must first construct enough central characters,
which we do by 
considering {\em dual Specht modules} following  
\cite[$\S$6]{schur}.
For a partition $\lambda$ of $d$, let $S^\lambda$ denote the usual Specht
module for the symmetric group $S_d$ over the field $F$ and
$S_\lambda = (S^\lambda)^*$ be its dual.
Given any $q\in F$, 
we can extend $S_\lambda$ to 
a module 
over the degenerate affine Hecke algebra $H_d$ so that $x_1$ acts by scalar multiplication by $q$. We denote the resulting $H_d$-module by $S_\lambda^q$.
If $d=d'+d''$, there is a natural embedding of $H_{d'} \otimes H_{d''}$ into 
$H_d$, so it makes sense to define the product
$$
M' \circ M'' = H_d \otimes_{H_{d'}\otimes H_{d''}} (M' \boxtimes M'')
$$
of an $H_{d'}$-module $M'$ and an $H_{d''}$-module
$M''$, where $\boxtimes$ denotes outer tensor product.
Given an $l$-multipartition
$\blambda = (\lambda^{(1)}, \lambda^{(2)},\dots,\lambda^{(l)})$
of $d$, the $H_d$-module
$$
S^{\bq}_{\blambda} := S^{q_1}_{\lambda^{(1)}}
\circ\cdots\circ S^{q_l}_{\lambda^{(l)}}
$$
factors through the quotient $H_d^{\bq}$ to give a well-defined
$H_d^{\bq}$-module. This is the {\em dual Specht module}
parametrized by the multipartition $\blambda$.

Let us 
compute the central 
character of the dual Specht module $S^{\blambda}_{\bq}$.
Note that if $M'$ is an $H_{d'}$-module of central character
$\chi(\bi')$ and 
$M''$ is an $H_{d''}$-module of central character $\chi(\bi'')$,
then 
$M' \circ M''$ is of central character
$\chi(\bi' \circ \bi'')$ where $\bi'\circ\bi''$ denotes 
the concatenation
$(i_1',\dots,i_{d'}', i_1'',\dots,i_{d''}')$.
This reduces to the problem of computing the central character
simply of $S^\lambda_q$ for a partition $\lambda$ of $d$, which is
well known:
for each $i,j \geq 1$
fill the box in the $i$th row and $j$th column 
of the Young diagram of $\lambda$
with the {\em residue} $(q + i-j)$;
then $S^q_\lambda$ is of central character 
parametrized by the tuple
$\bi^q_\lambda$
obtained by reading off
the 
residues in all the boxes
in some order. 
For example, if $\lambda = (4,2,1)$ and $q=5$ then the residues are
$$
\diagram{5&6&7&8\cr 
4&5\cr3\cr}
$$
so $S_\lambda^q$ is of central character
parametrized by $\bi_\lambda^q \sim (5,6,7,8,4,5,3)$.
Given $\blambda = (\lambda^{(1)},\dots,\lambda^{(l)})$,
we deduce that the central character of $S_{\blambda}^{\bq}$
is $\chi(\bi_\blambda^\bq)$ where
$$
\bi_{\blambda}^{\bq} = \bi_{\lambda^{(1)}}^{q_1} \circ \cdots\circ
\bi_{\lambda^{(l)}}^{q_l}.
$$
In this way, we have proved the existence of 
many 
central characters of $H_d^\bq$.

Now we proceed like in finite group theory.
Let $\bar R$ be a Noetherian domain
with maximal ideal $\bar{\mathfrak{m}}$
such that $F = \bar R / \bar{\mathfrak{m}}$
and the field of fractions of $\bar R$ is of characteristic $0$.
Let $\hat q_1,\dots,\hat q_l \in \bar R$ be lifts of the
parameters $q_1,\dots,q_l \in F$.
Let $R$ be the ring obtained by first localizing
the polynomial algebra
$\bar R[t_1,\dots,t_l]$ at the maximal ideal generated
by $\bar{\mathfrak m}$ and $t_1-\hat q_1,\dots,t_l-\hat q_l$,
and then completing with respect to the image of this maximal ideal.
We still have that $F = R  / \mathfrak{m}$, where
$\mathfrak{m}$ is the unique maximal ideal of $R$.
Also let $K$ be the field of fractions of $R$.
Letting $\bt = (t_1,\dots,t_l)$,
define
$H_d^\bt$ and $\mathcal H_d^\bt$
to be the degenerate cyclotomic Hecke algebras 
defined by the polynomial
$f(x) = (x-t_1)\cdots (x-t_l)$
over the
field $K$ and over
the ring $R$, respectively.
In view of the PBW theorem for degenerate cyclotomic Hecke
algebras, $H_d^{\bt}$ 
is naturally isomorphic to $K \otimes_{R}
\mathcal H_d^\bt$, and
$H_d^{\bq}$ is naturally isomorphic to 
$F \otimes_{R} \mathcal H_d^\bt$, viewing $F$ here as an
$R$-module so that each $t_i$ acts as multiplication by
$q_i$. 
Note moreover that the definition of dual Specht modules 
carries over unchanged to give modules
$S_{\blambda}^{\bt}$ for $H_d^{\bt}$ and $\mathcal S_{\blambda}^{\bt}$
for $\mathcal H_d^{\bt}$ for each $\blambda \in \MP_d(l)$, such that
$S_{\blambda}^{\bt} \cong 
K \otimes_{R} 
\mathcal S_{\blambda}^{\bt}$
and $S_{\blambda}^{\bq} \cong
F \otimes_{R} \mathcal S_{\blambda}^{\bt}$.
The following lemma is well known, but the proof given here is 
quite instructive.

\begin{Lemma}\label{easy}
The algebra $H_d^\bt$ is split semisimple. Moreover,
the dual Specht modules $S_{\blambda}^{\bt}$ give a complete set of
pairwise non-isomorphic irreducible $H_d^{\bt}$-modules.
\end{Lemma}

\begin{proof}
For each $i=1,\dots,l$,
let $H_d^{t_i}$ denote the 
degenerate cyclotomic Hecke algebra over $K$
defined by the polynomial $f(x) = (x-t_i)$.
There is an isomorphism $H_d^{t_i} 
\stackrel{\sim}{\rightarrow} K S_d$ 
which is the identity on $S_d$ and maps $x_1$ to $t_i$.
Since $K$ is a field of charcteristic zero,
we get from this that 
each $H_d^{t_i}$ is a split semisimple algebra and, by the classical 
representation theory of the symmetric group, the
dual Specht
modules $S_\lambda^{t_i}$ for all partitions $\lambda$ of $d$ give a 
complete set of pairwise non-isomorphic irreducible $H_d^{t_i}$-modules.
Since $t_1,\dots,t_l$ are algebraically independent,
the proof of \cite[Corollary 5.20]{schur} shows that
there is an isomorphism
$$
H_d^\bt \cong 
\bigoplus_{d_1+\cdots+d_l = d}
H_{d_1}^{t_1} \otimes\cdots\otimes H_{d_l}^{t_l}
$$
under which $S^{\bt}_{\blambda}$ corresponds to the outer tensor product
$S^{t_1}_{\lambda^{(1)}} \boxtimes\cdots\boxtimes
S^{t_l}_{\lambda^{(l)}}$ of dual Specht modules.
The lemma follows.
\end{proof}

Lemma~\ref{easy} implies that all the 
dual Specht modules $\{S^\bt_\blambda\mid\blambda \in X_d^\bt\}$
have different central characters.
One can also see this directly 
by observing from the combinatorial definition
that the tuples
$\bi^\bt_{\blambda}$ for $\blambda \in \MP_d(l)$ are
in different $\sim$-equivalence classes, i.e. in the generic case
the map
$$
\MP_d(l) \rightarrow X_d^\bt,
\qquad
\blambda \mapsto \bi^\bt_\blambda
$$
is injective. Actually, it is a bijection, by a trivial special
case of the 
following lemma completing the classification of blocks
of $H_d^\bq$ in general.

\begin{Lemma}\label{easier}
$X_d^\bq = \{\bi_{\blambda}^{\bq}\mid \blambda \in \MP_d(l)\}$.
\end{Lemma}

\begin{proof}
We have already noted that all $\bi_{\blambda}^{\bq}$ belong to 
$X_d^{\bq}$.
Conversely, we need show that the $\chi(\bi_{\blambda}^{\bq})$'s 
give all of the central characters of $H_d^{\bq}$.
This follows from the following claim:
we have that
$\prod_{\blambda \in \MP_d(l)} (z-\chi(\bi_{\blambda}^\bq)(z)) = 0$
for every $z \in Z(H_d^{\bq})$.
Note the following diagram commutes
$$
\begin{CD}
Z(\mathcal H_d^{\bt}) &@>\chi({\bi_{\blambda}^{\bt}})>>&R\\
@VVV&&@VVV\\
Z(H_d^{\bq}) &@>>\chi({\bi_{\blambda}^{\bq}})>&F
\end{CD}
$$
where the vertical maps are defined by evaluating
each $t_i$ at $q_i$.
So the claim follows if we can show that 
$\prod_{\blambda\in\MP_d(l)} (z-\chi(\bi_{\blambda})(z)) = 0$
for every 
$z \in \mathcal H_d^{\bt}$. But we have that $\mathcal H_d^{\bt} 
\subseteq H_d^{\bt}$, and
in the semisimple algebra 
$H_d^{\bt}$ it is certainly the case 
that
$\prod_{\blambda\in\MP_d(l)} (z-\chi(\bi_{\blambda}^{\bt})(z)) = 0$
because the $\chi({\bi_{\blambda}^\bt})$'s
for all $\blambda \in \MP_d(l)$ are the central characters
of a full set of irreducible $H_d^{\bt}$-modules, thanks to Lemma~\ref{easy}.
\end{proof}

For $\bi \in X_d^\bq$, let $b(\bi)$ be the
primitive central idempotent corresponding to the central
character $\chi(\bi)$,
that is, $b(\bi)$ is the unique element of $Z(H_d^\bq)$ that 
acts as one on irreducible modules of central
character $\chi(\bi)$ and as zero on all other irreducibles.
Thus, we have that
$$
H_d^\bq = \bigoplus_{\bi \in X_d^{\bq}} b(\bi) H_d^\bq.
$$
This is the decomposition of $H_d^\bq$ into blocks.
Similarly, recalling the bijection
$\MP_d(l) \rightarrow X_d^\bt, \blambda \mapsto \bi_\blambda^\bt$,
we can define idempotents $b(\blambda) \in Z(H_d^\bt)$
for each $\blambda \in \MP_d(l)$
such that $b(\blambda)$ acts as
one on $S^\bt_\blambda$ 
and as zero on all other dual Specht modules.
Of course, the resulting decomposition
$$
H_d^\bt = \bigoplus_{\blambda \in \MP_d(l)} b(\blambda) H_d^\bt
$$
is the Wedderburn decomposition of the semisimple algebra
$H_d^\bt$.

Since $R$ is a Noetherian ring complete
with respect to the maximal ideal $\mathfrak{m}$,
and moreover we know that 
$Z(\mathcal H_d^\bt)$ surjects onto $Z(H_d^\bq)$ by Theorem 1,
there is a unique lift of each $b(\bi) \in Z(H_d^\bq)$ to a central idempotent
$\hat b(\bi) \in \mathcal H_d^\bt$; see e.g. \cite[Corollary 7.5]{E}.
This lifts the block decomposition of $H_d^\bq$
to a decomposition
$$
H_d^\bt = \bigoplus_{\bi \in X_d^{\bq}}
\hat b(\bi) H_d^\bt
$$
of the semisimple algebra $H_d^\bt$.
Finally, the commutative diagram from the proof
of Lemma~\ref{easier} implies for each $\bi \in X_d^\bq$ that
$$
\hat b(\bi) H_d^\bt = 
\bigoplus_{\stackrel{\blambda \in \MP_d(l)}{\bi^{\bq}_{\blambda} = \bi}}
b(\blambda) H_d^\bt.
$$
Now we can prove the only new result of the section, as follows.

\begin{Theorem}\label{easiest}
For $\bi \in X_d^{\bq}$, 
the dimension of the center of the block 
$b(\bi) H_d^{\bq}$ is equal to the 
number of $l$-multipartitions
$\blambda$ of $d$ 
such that $\bi^{\bq}_{\blambda} = \bi$.
\end{Theorem}

\begin{proof}
By Theorem 1, $Z(\mathcal H_d^\bt)$ is a free $R$-module of finite rank.
Since $R$ is a local ring, it follows that the summand
$Z(\hat b(\bi) \mathcal H_d^\bt) = 
\hat b(\bi) Z(\mathcal H_d^\bt)$ is also free,
of rank equal to
$\dim_{K} Z(\hat b(\bi) H_d^{\bt})$.
Since each $b(\blambda) H_d^\bt$ 
is a full matrix algebra with a one dimensional center, 
we know from the preceeding discussion that
this dimension is equal to the number of $l$-multipartitions
$\blambda$ with $\bi^\bq_\blambda = \bi$.
By Theorem 1 again, the 
isomorphism $F \otimes_R \mathcal H_d^\bt \stackrel{\sim}{\rightarrow}
H_d^\bq$
induces an isomorphism
$F \otimes_R Z(\mathcal H_d^\bt)
\stackrel{\sim}{\rightarrow} Z(H_d^\bq)$.
From this, we get an isomorphism
$F \otimes_R Z(\hat b(\bi) \mathcal H_d^\bt)
\stackrel{\sim}{\rightarrow} 
Z(b(\bi) H_d^\bq)$.
So $\dim_F Z(b(\bi) H_d^\bq)$ is the same as the
rank of $Z(\hat b(\bi) \mathcal H_d^\bt)$,
i.e. the number of $l$-multipartitions $\blambda$ of $d$
with $\bi^\bq_\blambda = \bi$.
\end{proof}

\section{The center of parabolic category $\mathcal O$}

Let $\mathfrak{g} = \mathfrak{gl}_n(\C)$ with natural module $V$.
We denote the standard basis for $V$ by $v_1,\dots,v_n$ and 
use the notation $e_{i,j}$ for the matrix units in 
$\mathfrak{g}$.
Let $\mathfrak{d}$ be the subalgebra of $\mathfrak{g}$ of diagonal matrices
and $\mathfrak{b}$ be the standard Borel subalgebra of upper triangular 
matrices.
Let $\eps_1,\dots,\eps_n$ be the basis for $\mathfrak{d}^*$
dual to the standard basis $e_{1,1},\dots,e_{n,n}$ for $\mathfrak{d}$.
We write $L(\alpha)$ for the irreducible highest weight module
of highest weight $(\alpha-\rho)$, where
$\rho$ is the weight $-\eps_2 - 2 \eps_3 - \cdots - (n-1) \eps_n$.
Viewing elements of $S(\mathfrak{d})$ as polynomial functions on $\mathfrak{d}^*$,
the {\em Harish-Chandra homomorphism}
$$
\Psi:Z(\mathfrak{g}) \stackrel{\sim}{\longrightarrow} S(\mathfrak{d})^{S_n}
$$
can be defined by declaring that $\Psi(z)$ is the unique element
of $S(\mathfrak{d})$ with the property
that $z$ acts on $L(\alpha)$ by the scalar $(\Psi(z))(\alpha)$ for
each $\alpha \in \mathfrak{d}^*$.
Its image is the algebra $S(\mathfrak{d})^{S_n}$ of symmetric polynomials (for 
the usual action of $S_n$ on $\mathfrak{d}$ not the dot action).

Letting $\tilde e_{i,j} := e_{i,j} + \delta_{i,j}(u+1-i)$ for short,
it is classical that the coefficients $z_1,\dots,z_n$ of the polynomial
$$
z(u) = \sum_{r=0}^n z_r u^{n-r}
:=
\sum_{w \in S_n} \operatorname{sgn}(w) \tilde e_{w1,1} \cdots
\tilde e_{wn, n}
\in U(\mathfrak{g})[u]
$$
are algebraically independent generators
for the center $Z(\mathfrak{g})$ of $U(\mathfrak{g})$.
We adopt the convention that $z_r = 0$ for $r > n$.
The image of $z(u)$ under $\Psi$ is given by the formula
$$
\Psi(z(u)) = (u+e_{1,1}) \cdots (u+e_{n,n}).
$$
Hence, for $\alpha = \sum_{i=1}^n a_i \eps_i \in \mathfrak{d}^*$,
the central element $z_r$ acts on $L(\alpha)$ as the scalar
$e_r(\alpha) = e_r(a_1,\dots,a_n)$, 
the $r$th {\em elementary symmetric function}
evaluated at the numbers $a_1,\dots,a_n$.
Let $P$ denote the free abelian group on basis $\{\gamma_a\mid a
\in \C\}$.
For $\alpha = \sum_{i=1}^n a_i \eps_i \in \mathfrak{d}^*$,
let $\nu(\alpha)=\gamma_{a_1}+\cdots+\gamma_{a_n} \in P$.
The point of this definition is that
$L(\alpha)$ and $L(\beta)$ have the same central character
if and only if $\nu(\alpha) = \nu(\beta)$.
In this way, we have parametrized the central characters
of $U(\mathfrak{g})$ by the set of all
$\nu \in P$ whose coefficients are non-negative integers summing to
$d$.

Let $\Delta:U(\mathfrak{g}) \rightarrow U(\mathfrak{g}) \otimes U(\mathfrak{g})$ be the canonical comultiplication on the universal enveloping algebra of
$\mathfrak{g}$.
We are only going to need to work with
the homomorphism $\delta:U(\mathfrak{g}) \rightarrow U(\mathfrak{g}) \otimes 
\End_\C(V)$ obtained by composing $\Delta$ with the map $1 \otimes \phi$
where $\phi:U(\mathfrak{g}) \rightarrow \End_\C(V)$ here 
is the algebra homomorphism
arising from the representation of $\mathfrak{g}$ on $V$.
Also, let
$$
\Omega = \sum_{i,j=1}^n e_{i,j} \otimes e_{j,i}
\:\in U(\mathfrak{g}) \otimes \End_{\C}(V).
$$
The following lemma is probably classical.

\begin{Lemma}\label{ought}
For $r \geq 0$, we have that
$$
\delta(z_r) = z_r\otimes 1
+ \sum_{s=0}^{r-1} (-1)^s (z_{r-1-s} \otimes 1) \Omega^s.
$$
\end{Lemma}

\begin{proof}
Both sides of the equation are elements of
$U(\mathfrak{g}) \otimes \End_\C(V)$, so can be viewed as $n \times n$ matrices
with entries in $U(\mathfrak{g})$.
To see that these matrices are equal, it suffices to check that 
their entries act in the same way on
sufficiently many finite dimensional 
irreducible representations of $\mathfrak{g}$.
This reduces to the following problem.
Take $\alpha = \sum_{i=1}^n a_i \eps_i \in \mathfrak{d}^*$ such that $L(\alpha)$ 
is finite 
dimensional and 
$$
L(\alpha) \otimes V \cong \bigoplus_{i=1}^n
L(\alpha + \eps_i).
$$
We need to show that the left and right hand sides of the given equation
define the same endomorphism of $L(\alpha) \otimes V$.
For such an $\alpha$, let $M := L(\alpha) \otimes V$ and write
$v_+$ for a highest weight vector in $L(\alpha)$.
For $i=0,\dots,n$, define $M_i$ to be the submodule of $M$
generated by the vectors $v_+ \otimes v_j\:(j \leq i)$.
Since $v_+ \otimes v_i$ is a highest weight vector 
of weight $\alpha + \eps_i$ modulo $M_{i-1}$,
the assumption on $\alpha$ implies that
$0 = M_0 \subset M_1 \subset \cdots \subset M_n = M$ 
is a filtration of $M$
such that $M_i / M_{i-1} \cong L({\alpha+\eps_i})$.
Since the filtration splits, there is a unique highest weight
vector $x_i \in M$ such that $x_i \equiv v_+ \otimes v_i \pmod{M_{i-1}}$.
Now we just check that the left and right hand sides of the given 
equation act on these highest weight
vectors by the same scalar for each $i=1,\dots,n$.
Of course $\delta(z_r)$ acts on $x_i$ as $e_r(\alpha+\eps_i)$,
while each $z_t \otimes 1$ acts as $e_t(\alpha)$ on all of $M$.
Finally, since $\Omega$ defines a $\mathfrak{g}$-module
endomorphism of $M$, it leaves $M_{i-1}$ invariant and maps
$x_i$ to a scalar multiple of itself.
To compute the scalar, note that
\begin{align*}
\Omega (v_+ \otimes v_i)
&= \sum_{j\leq i} (e_{i,j} v_+) \otimes v_j=
(e_{i,i}v_+) \otimes v_i + 
\sum_{j < i} (e_{i,j} v_+) \otimes v_j
\\
&= (a_i+i-1) v_+ \otimes v_i
+ \sum_{j < i} (e_{i,j} (v_+ \otimes v_j)
- v_+ \otimes v_i)\\
&\equiv a_i v_+ \otimes v_i \pmod{M_{i-1}}.
\end{align*}
Hence, $\Omega x_i = a_i x_i$.
So the equation we are trying to prove reduces to checking that
$$
e_r(\alpha+\eps_i) = e_r(\alpha) + \sum_{s=0}^{r-1} (-1)^s e_{r-1-s}(\alpha) a_i^s
$$
for each $i=1,\dots,n$.
This follows from the following general identity
which is true for all $r,k \geq 0$:
$$
e_r(u_1,\dots,u_k,u+1) 
=
e_r(u_1,\dots,u_k,u)+\sum_{s=0}^{r-1} (-1)^s e_{r-1-s}(u_1,\dots,u_k,u) u^s.
$$
To see this, expand both sides
using the obvious formula
$e_t(u_1,\dots,u_k, v) = e_t(u_1,\dots,u_k)
 + e_{t-1}(u_1,\dots,u_k) v$. 
\iffalse
note that
$$
e_{r-1}(u_1,\dots,u_k) = \sum_{s=0}^{r-1} (-1)^s e_{r-1-s} (u_1,\dots,u_k,u)
u^s,
$$
as follows by expanding the right hand side
using 
Hence,
\begin{align*}
e_r(u_1,\dots,u_k,u+1) 
&=
e_r(u_1,\dots,u_k) + e_{r-1}(u_1,\dots,u_k) u + e_{r-1}(u_1,\dots,u_k)\\
&=
e_r(u_1,\dots,u_k, u) + e_{r-1}(u_1,\dots,u_k)\\
&=
e_r(u_1,\dots,u_k,u)+\sum_{s=0}^{r-1} (-1)^s e_{r-1-s}(u_1,\dots,u_k,u) u^s,
\end{align*}
as required.
\fi
\end{proof}

Let $M$ be any $\mathfrak{g}$-module.
Recall from \cite[$\S$2.2]{AS} that
the degenerate affine Hecke algebra
$H_d$ over the ground field $\C$ acts naturally on the right on $M \otimes V^{\otimes d}$
by $\mathfrak{g}$-module endomorphisms.
The action of each $w \in S_d$ arises from its usual
action on $V^{\otimes d}$ by place permutation.
The action of $x_1$ (from which one can deduce the action of
all other $x_i$'s) is the same as left multiplication by
$\Omega \otimes 1^{\otimes (d-1)}$.
For any partition
$\mu$ with $\ell(\mu) \leq d$,
recall the notation $p_d(\mu)$ introduced just before Theorem~\ref{pr};
we are now viewing this expression as an element of $H_d$.

\begin{Lemma}\label{rd}
For any $r,d \geq 0$
and any highest weight module $M$ of highest weight
$\alpha-\rho \in \mathfrak{d}^*$, the endomorphism of
$M \otimes V^{\otimes d}$ define by
left multiplication by $z_r \in Z(\mathfrak{g})$
is equal to the endomorphism defined by right multiplication by
$$
\sum_{\mu} (-1)^{|\bar\mu|}\binom{d-\ell(\bar\mu)}{d-\ell(\mu)}
e_{r-|\mu|}(\alpha) p_d(\bar\mu)\:\in Z(H_d)
$$
where the sum is over partitions $\mu$ of length $\ell(\mu) \leq d$
and size $|\mu| \leq r$,
and $\bar\mu$ denotes the partition 
$(\mu_1-1 \geq \dots\geq \mu_{\ell(\mu)}-1)$
obtained from $\mu$ by removing the first
column of its diagram.
\end{Lemma}

\begin{proof}
Let
$\delta_d: U(\mathfrak{g}) \rightarrow 
U(\mathfrak{g}) \otimes \End_\C(V)^{\otimes d}$
be the map defined inductively by setting $\delta_0 = 1$ and
$\delta_d = (\delta \otimes 1^{\otimes (d-1)}) \circ 
\delta_{d-1}$ for $d \geq 1$.
Let $\Omega_i := (\delta_{i-1} \otimes 1)(\Omega) \otimes 1^{(d-i)}
\in 
U(\mathfrak{g}) \otimes \End_\C(V)^{\otimes d}$.
If we adopt the convention that
$(-\Omega)^{-1} = 1$, 
we can write the conclusion of Lemma~\ref{ought} simply as
$$
\delta(z_r) = \sum_{s=0}^r (z_{r-s} \otimes 1) (-\Omega)^{s-1}.
$$
Proceeding from this by induction on $d$, 
it is straightforward to deduce that
$$
\delta_d(z_r) = 
\sum_{\substack{s_1,\dots,s_d \geq 0 \\ s_1+\cdots+s_d \leq r}}
(z_{r-s_1-\cdots-s_d} \otimes 1^{\otimes d}) 
 (-\Omega_1)^{s_1-1} \cdots (-\Omega_d)^{s_d-1},
$$
interpreting the right hand side with same convention.
Since $x_{i+1} = s_i x_i s_i+s_i$
and $x_1$ acts as $\Omega_1$ by definition,
one checks by induction that $x_i$ acts as $\Omega_{i}$ for each $i$.
Hence on applying our expression 
to $M \otimes V^{\otimes d}$, we deduce that
$z_r$ acts in the same way as 
$$
\sum_{\substack{s_1,\dots,s_d \geq 0 \\ s_1+\cdots+s_d \leq r}}
e_{r-s_1-\cdots-s_d}(\alpha) (-x_1)^{s_1-1} \cdots (-x_d)^{s_d-1},
$$
again interpreting $(-x_i)^{-1}$ as $1$.
It is now a combinatorial exercise to rewrite this expression as formulated
in the statement
of the lemma.
\end{proof}

\begin{Corollary}\label{dwarf}
For any highest weight module $M$,
the subalgebra of the algebra
$\End_{\C}(M \otimes V^{\otimes d})^{\op}$
generated by the endomorphisms $z_r\:(1 \leq r \leq n)$ 
coincides with the subalgebra generated by 
$x_1^r+\cdots+x_d^r\:(1 \leq r \leq d)$.
\end{Corollary}

\begin{proof}
Since we are working over a field of characteristic $0$,
any symmetric polynomial in $x_1,\dots,x_d$
lies in the subalgebra generated by the power sums
$x_1^r + \cdots + x_d^r\:(1 \leq r \leq d)$.
By Lemma~\ref{rd}, the endomorphism defined by $z_r$ can be expressed as
a symmetric polynomial in $x_1,\dots,x_d$, 
so it lies in the subalgebra generated by the power sums.

Conversely, we show by induction on $r \geq 0$ that every homogeneous
symmetric polynomial in $x_1,\dots,x_d$ of degree $r$ 
acts on $M \otimes V^{\otimes d}$ in the same way
as some element of the subalgebra generated by
$z_1,\dots,z_{r+1}$.
For the induction step, every homogeneous symmetric
polynomial of degree $r$ lies in the subalgebra generated
by the power sums $x_1^s+\cdots+x_d^s\:(1 \leq s \leq r)$.
By induction all of these power sums with $s < r$ certainly 
lie in the subalgebra
generated by $z_1,\dots,z_{r+1}$,
so it just remains to show that $x_1^r+\cdots+x_d^r$ does too.
By Lemma~\ref{rd}, the image of $z_{r+1}$ is the same as the image of
$x_1^r + \cdots + x_d^r$ (which is the term $p_d(\bar\mu)$ when $\mu = (r+1)$) 
plus a linear combination 
of symmetric polynomials in $x_1,\dots,x_d$ of strictly smaller degree,
which we already have by the induction hypothesis. 
\end{proof}

Let $\mu = (\mu_1,\dots,\mu_l)$ be a composition of $n$ and
let $\mathfrak{p}$ be the corresponding standard parabolic subalgebra of $\mathfrak{g}$ with standard Levi subalgebra $\mathfrak{h}
\cong \mathfrak{gl}_{\mu_1}(\C)\oplus\cdots\oplus \mathfrak{gl}_{\mu_l}(\C)$,
as in the introduction.
We are interested in the category $\mathcal O^\mu$
of all finitely generated
$\mathfrak{g}$-modules that are locally finite over $\mathfrak{p}$
and integrable over $\mathfrak{h}$.
Also let $\bq = (\mu_1,\dots,\mu_l)$ and let $H_d^{\bq}$ denote
the degenerate cyclotomic Hecke algebra
from the previous section
over the ground field $F = \C$. 
We are going to apply
the Schur-Weyl duality for higher levels from \cite{schur} 
(taking the choice of origin there to be $\bc = (n,\dots,n)$)
to connect the category $\mathcal O^\mu$
to the finite dimensional 
algebras $H_d^\bq$ for all $d \geq 0$.
Actually, \cite{schur} only considered the special case that
$\mu$ is a partition, 
i.e. $\mu_1 \geq \cdots \geq \mu_l$, so we will need to
 extend some of the results of \cite{schur} to the general case
as we go.

To start with, we need some combinatorial definitions. Let
$$
\Col^\mu = 
\left\{\alpha = \sum_{i=1}^n a_i \eps_i \in \mathfrak{d}^*
\:\bigg|\:
\begin{array}{ll}
a_1,\dots,a_n \in \Z \text{ such that }
a_i > a_{i+1}\text{ for}\\\text{all }i\neq \mu_1,\mu_1+\mu_2,\dots,\mu_1+\cdots+\mu_l\end{array}\right\},
$$
so called because its elements can be visualized as column strict tableaux
of column shape $\mu$ like in \cite[$\S$2]{dual}.
The irreducible modules in $\mathcal O^\mu$
are the modules $\{L(\alpha)\mid \alpha \in \Col^\mu\}$.
Hence the set $Y^\mu = \{\nu(\alpha)\mid \alpha \in \Col^\mu\}$ 
naturally parametrizes the central characters arising from
modules in $\mathcal O^\mu$.
Given $\nu \in Y^\mu$, 
we let $\Col^\mu_\nu 
= \{\alpha \in \Col^\mu\mid \nu(\alpha) = \nu\}$
and define
$\mathcal O^\mu_\nu$ to be
the Serre subcategory of $\mathcal O^\mu$
generated by the modules 
$\{L(\alpha)\mid \alpha \in \Col^\mu_\nu\}$.
The category $\mathcal O^\mu$ then decomposes as
$$
\mathcal O^\mu = \bigoplus_{\nu \in Y^\mu} \mathcal O^\mu_\nu.
$$
This is the same as the central character decomposition of $\mathcal O^\mu$
that was described in the introduction.

Let $\gamma \in Y^\mu$ be the special element
$\gamma = \sum_{i=1}^l \sum_{a=1}^{\mu_i} \gamma_a$.
The key feature of $\gamma$ is that the set
$\Col^\mu_\gamma$ contains just one weight $\alpha$.
In other words, for this $\alpha$, 
$L(\alpha)$ is the unique
irreducible module in $\mathcal O^\mu$
with central character parametrized by $\gamma$.
This special irreducible module, which we denote henceforth by $P^\mu$,
is automatically projective.
\iffalse

weight $\alpha \in \Col^\mu$ such that 
$\nu(\alpha) = \gamma$, namely,
the weight
$$
\alpha = \sum_{i=1}^l \sum_{a=1}^{\mu_i}
a \eps_{\mu_1+\cdots+\mu_i+1-a}.
$$
For this $\alpha$, the module $L(\alpha)$ 
is both irreducible and projective.
We denote it by $P^\mu$ from now on.
\fi
For $d \geq 0$, 
let 
\begin{align*}
Y^\mu_d &= \{\gamma - (\gamma_{i_1}-\gamma_{i_1+1}) - \cdots
- (\gamma_{i_d}-\gamma_{i_d+1}) \in Y^\mu
\mid i_1,\dots,i_d \in \Z, \},\\
\Col^\mu_d &= \{\alpha \in \Col^\mu\mid \nu(\alpha) \in Y^\mu_d\}.
\end{align*}
\iffalse
,\\
\mathcal O^\mu_d &= \bigoplus_{\nu \in \Col^\mu_d} \mathcal O^\mu_\nu.
\fi
The irreducible modules $\{L(\alpha)\mid \alpha \in \Col^\mu_d\}$
are significant because they are 
exactly the irreducible constituents of the
module $P^\mu \otimes V^{\otimes d}$.
This statement is proved in \cite[$\S$4]{schur}
in the case that $\mu$ is a partition, 
and the same argument works in general.

\begin{Lemma}\label{comb}
The map $Y^\mu_d \rightarrow X^{\bq}_d$ sending
$\nu = \gamma - (\gamma_{i_1}-\gamma_{i_1+1})
- \cdots - (\gamma_{i_d}-\gamma_{i_d+1})$
to $\bi = (i_1,\dots,i_d)$
is injective with image equal to
$$
\{\bi^\bq_\blambda\mid \blambda \in \MP_d(l)
\text{ such that }\ell(\lambda^{(r)}) \leq \mu_r
\text{ for $r = 1,\dots,l$}\}.
$$
For $\nu \in Y^\mu_d$ 
corresponding to $\bi \in X^{\bq}_d$ in this way,
the map
$\Col^\mu_\nu \rightarrow \MP_d(l)$
sending $\alpha = \sum_{i=1}^n a_i \eps_i$ to the multipartition
$\blambda = (\lambda^{(1)},\dots,\lambda^{(l)})$ such that
$\lambda^{(r)} = (a_{\mu_1+\cdots+\mu_{r-1}+1} - \mu_r,
a_{\mu_1+\cdots+\mu_{r-1}+2} - \mu_{r-1},
\dots,
a_{\mu_1+\cdots+\mu_r} - 1)$
is injective with image equal to
$\{\blambda \in \MP_d(l)\mid \bi^\bq_\blambda  =\bi\}$.
\end{Lemma}

\begin{proof}
We leave this as a simply combinatorial exercise.
It is helpful to use the interpretation of $\Col^\mu_\nu$
as the set of column strict tableaux of column shape $\mu$
and type $\nu$ as in \cite[$\S$4]{schur}.
\end{proof}

As explained before Lemma~\ref{rd}, 
the degenerate affine Hecke algebra $H_d$ acts 
on the right on $P^\mu \otimes V^{\otimes d}$
by $\mathfrak{g}$-module endomorphisms.
We let
$$
\rho^\mu:H_d \rightarrow \End_{\mathfrak{g}}(P^\mu \otimes V^{\otimes d})^{\operatorname{op}}
$$
be the resulting algebra homomorphism.
The following lemma is the key to extending the results from \cite{schur}
to general $\mu$.

\begin{Lemma}\label{ms}
Given another composition $\mu' \sim \mu$,
there is an algebra isomorphism $\iota_{\mu,\mu'}$ 
making the following diagram commute:
$$
\begin{CD}
&\: H_d &\\
\\
\End_{\mathfrak{g}}(P^\mu \otimes V^{\otimes d})^{\op}
@>\sim > \iota_{\mu,\mu'}> &\End_{\mathfrak{g}}(P^{\mu'} \otimes V^{\otimes d})^{\op}.
\end{CD}
\begin{picture}(0,0)
\put(-149,1){\makebox(0,0){$\swarrow$}}
\put(-153,11){\makebox(0,0){$_{\rho^\mu}$}}
\put(-149,1){\line(1,1){18}}
\put(-80,1){\makebox(0,0){$\searrow$}}
\put(-74,11){\makebox(0,0){$_{\rho^{\mu'}}$}}
\put(-80,1){\line(-1,1){18}}
\end{picture}
$$
Moreover, $\iota_{\mu,\mu'}$ intertwines the natural actions
of $Z(\mathfrak{g})$ on the two endomorphism algebras.
\end{Lemma}

\begin{proof}
We appeal to an argument due to Mazorchuk and Stroppel.
By the proof of \cite[Theorem 5.4]{MS},
there is an adjoint pair $(F,G)$ of functors between the bounded derived categories
$$
D^b(\mathcal O^\mu)\quad\qquad
D^b(\mathcal O^{\mu'}).
\begin{picture}(0,15)
  \put(-59,4){\makebox(0,0){$\stackrel{\stackrel{\,\scriptstyle{F}_{\phantom{,}}}{\displaystyle\longrightarrow}}{\stackrel{\displaystyle\longleftarrow}{\scriptstyle G}}$}}
\end{picture}
$$
with the following properties:
\begin{itemize}
\item[(i)] $F$ and $G$ commute with
tensoring with finite
dimensional $\mathfrak{g}$-modules, that is, for any finite dimensional
$\mathfrak{g}$-module $V$ there are given
natural isomorphisms $\alpha_V: F \circ 
? \otimes V {\rightarrow}
? \otimes V \circ F$ and $\beta_V:G \circ ? \otimes V 
{\rightarrow} ? \otimes V \circ G$ such that the following two diagrams commute
for any morphism $f:V \rightarrow W$ of finite dimensional 
$\mathfrak{g}$-modules:
$$
\begin{CD}
F \circ ? \otimes V &@>\alpha_V >>& ? \otimes V \circ F\\
@V 1 (\operatorname{id}\otimes f) VV&&@VV (\operatorname{id}\otimes f)1 V\\
F \circ ? \otimes W &@>>\alpha_W>&? \otimes W \circ F
\end{CD}
\qquad\qquad\begin{CD}
G \circ ? \otimes V &@>\beta_V >>& ? \otimes V \circ G\\
@V 1 (\operatorname{id}\otimes f) VV&&@VV (\operatorname{id}\otimes f)1 V\\
G \circ ? \otimes W &@>>\beta_W>&? \otimes W \circ G
\end{CD}
$$
\item[(ii)] The isomorphisms $\alpha_V$ and $\beta_V$ 
are compatible with the unit
$\eta:\operatorname{Id} \rightarrow G \circ F$ and counit
$\eps:F \circ G \rightarrow \operatorname{Id}$ 
of the canonical adjunction between $F$ and $G$, i.e. the following diagrams
commute:
$$
\begin{CD}
? \otimes V &@>1 \eta>>& ? \otimes V \circ G \circ F\\
@VV\eta 1 V&&@AA{\beta_V} 1 A\\
G \circ F \circ ? \otimes V &@>>1 \alpha_V>&G \circ ? \otimes V \circ F
\end{CD}
\quad
\begin{CD}
? \otimes V &@<1 \eps<<& ? \otimes V \circ F \circ G\\
@AA\eps 1 A&&@AA{\alpha_V1} A\\
F \circ G \circ ? \otimes V &@>>1 \beta_V>&F \circ ? \otimes V \circ G
\end{CD}
$$
\item[(iii)] $F$ and $G$ restrict to mutually inverse equivalences
of categories between $\mathcal O^\mu_\gamma$ and $\mathcal O^{\mu'}_\gamma$.
\item[(iv)] The following associativity pentagon commutes for any 
two finite dimensional $\mathfrak{g}$-modules $V$ and $W$:

\vspace{0.5mm}
$$
\begin{CD}
F \circ ? \otimes (V \otimes W) &&&&@>>>? \otimes(V \otimes W) \circ F\\
@|&&&&@|\\
F \circ ?\otimes W \circ ? \otimes V&@>>\alpha_W 1> ? \otimes W \circ F \circ ? \otimes V
&@>>1 \alpha_V > ? \otimes W \circ ? \otimes V \circ F
\end{CD}
\begin{picture}(0,0)
\put(-240,27.28){\line(1,0){128}}
\put(-160,35){\makebox(0,0){$_{\alpha_{V \otimes W}}$}}
\end{picture}
$$
\item[(v)] $F$ transforms the endomorphism 
of an object $M$
defined by left multiplication
by $z \in Z(\mathfrak{g})$ to the 
endomorphism of $FM$ defined by left multiplication by the same element $z$.
\end{itemize}

By (iii), we can choose an isomorphism $F(P^\mu) \cong P^{\mu'}$ allowing us to simply
identify $F(P^\mu)$ with $P^{\mu'}$.
The isomorphism $\alpha_{V^{\otimes d}}$ then allows us to 
identify $F(P^\mu \otimes V^{\otimes d})$ with $P^{\mu'} \otimes V^{\otimes d}$.
Applying $G$ to this, we get an isomorphism $G(F(P^\mu \otimes V^{\otimes d}) 
\cong G(P^{\mu'} \otimes V^{\otimes d})$,
hence on composing with the counit of the adjunction we 
get a map 
$P^\mu \otimes V^{\otimes d} \rightarrow G(P^{\mu'} \otimes V^{\otimes d})$ 
which by (ii) and (iii) is an isomorphism.
Using this, we also identify
$G(P^{\mu'} \otimes V^{\otimes d})$ with $P^\mu \otimes V^{\otimes d}$.
Now the functor $F$ defines an algebra homomorphism
$$
\iota_{\mu,{\mu'}}:\End_{\mathfrak{g}}(P^\mu \otimes V^{\otimes d})^{\op} 
\rightarrow \End_{\mathfrak{g}}(P^{\mu'}
\otimes V^{\otimes d})^{\op}
$$
and the functor $G$ defines a
homomorphism $$
\iota_{\mu',\mu}:\End_{\mathfrak{g}}(P^{\mu'}\otimes V^{\otimes d})^{\op} \rightarrow \End_{\mathfrak{g}}(P^\mu \otimes V^{\otimes d})^{\op}
$$
such that $\iota_{\mu,\mu'}$ and $\iota_{\mu',\mu}$ are mutual inverses.
Moreover, by (v), both homomorphisms intertwine the natural actions of
$Z(\mathfrak{g})$.

It just remains to check that $\iota_{\mu,\mu'}$ is compatible with the
action of $H_d$. The compatibility of $\iota_{\mu,\mu'}$ 
with 
the action of each $w \in S_d$ follows immediately from the naturality
in (i).
So it suffices to show that
it is compatible with the action of $x_1$. For this, we first reduce
using (iv) to checking compatiblity just in the special case $d=1$.
In that case it follows from (v) since by Lemma~\ref{ought}
we have that $x_1$ acts as left multiplication by
$\Omega = z_2 \otimes 1 + z_1 \otimes 1 - \delta(z_2)$,
and $z_2 \otimes 1$ and $z_1 \otimes 1$ act by
the same scalars on $P^\mu \otimes V$ and $P^{\mu'}\otimes V$.
\end{proof}

We can formulate the critical result needed from \cite{schur}
as follows.

\begin{Theorem}\label{dcp}
The image of
$\rho^\mu: H_d \rightarrow 
\End_{\C}(P^\mu \otimes V^{\otimes d})^{\op}$ coincides with the endomorphism algebra
$\End_{\mathfrak{g}}(P^\mu \otimes V^{\otimes d})^{\op}$.
Moreover, 
the representation 
$\rho^\mu$ factors through the quotient $H_d^\bq$ of $H_d$,
and the kernel of the induced map 
$H_d^\bq \twoheadrightarrow \End_{\mathfrak{g}}(P^\mu \otimes V^{\otimes d})^{\op}$
is generated by $(1-e)$, where $e \in H_d^\bq$ is the central idempotent
$e = \sum_{\bi} b(\bi)$
summing over all $\bi$ lying in the set
$$
\{\bi^{\bq}_{\blambda}\mid \blambda \in \MP_d(l)
\text{ such that }\ell(\lambda^{(r)}) \leq \mu_r \text{ for }
r=1,\dots,l\}
$$
from Lemma~\ref{comb}.
Hence, $\rho^\mu$ induces an isomorphism between
the sum of blocks
$e H_d^\bq$  of $H_d^\bq$ and the endomorphism algebra
$\End_{\mathfrak{g}}(P^\mu \otimes V^{\otimes d})^{\op}$.
\end{Theorem}

\begin{proof}
If $\mu$ is a partition, this follows by
\cite[Theorem 5.13]{schur}
and \cite[Corollary 6.7]{schur}.
It then follows for arbitrary $\mu$ too by
Lemma~\ref{ms}.
\end{proof}

\begin{Corollary}\label{C2}
The center of $\End_{\mathfrak{g}}(P^\mu \otimes V^{\otimes d})^{\op}$
is generated by the endomorphisms $z_1,\dots,z_n$.
\end{Corollary}

\begin{proof}
By Corollary~\ref{ps} we know already that $Z(H_d^\bq)$
is generated by the power sums 
$x_1^r+\cdots+x_d^r$ for $1 \leq r \leq d$.
By Theorem~\ref{dcp}, we can identify the endomorphism algebra
$\End_{\mathfrak{g}}(P^\mu \otimes V^{\otimes d})^{\op}$ with
$e H_d^\bq$ for some central idempotent $e \in H_d^\bq$.
So its center is generated
by the restrictions of these power sums to 
the module $P^\mu \otimes V^{\otimes d}$.
Now apply Corollary~\ref{dwarf}.
\end{proof}

Now let us restrict attention to a single block.
For the rest of the article, we fix $\nu \in Y^\mu_d$
and define 
$\bi = (i_1,\dots,i_d) \in X_d^\bq$
from $\nu = \gamma - (\gamma_{i_1}-\gamma_{i_1+1})-
\cdots - (\gamma_{i_d}-\gamma_{i_d+1})$ as in Lemma~\ref{comb}.
Let $e^\mu_\nu \in 
\End_{\mathfrak{g}}(P^\mu \otimes V^{\otimes d})^{\op}$ 
be the central
idempotent projecting $P^\mu \otimes V^{\otimes d}$
onto its component
of generalized central character
parametrized by $\nu$.
Note $(P^\mu \otimes V^{\otimes d}) e^\mu_\nu$ is non-zero;
see e.g. \cite[Lemma 4.2]{schur}.
Identifying $\End_{\mathfrak{g}}(P^\mu \otimes V^{\otimes d})^{\op}$
with $e H_d^\bq$ according to Theorem~\ref{dcp},
it follows that $e^\mu_\nu$ is identified with a non-zero 
central idempotent in $H_d^\bq$.

\begin{Lemma}
$e^\mu_\nu = b(\bi)$.
\end{Lemma}

\begin{proof}
We first prove this in the special case that $\mu$ is a partition.
Certainly $e^\mu_\nu$ 
is a non-zero sum of the primitive central idempotents
$b(\bi)$ for $\bi \in X_d^\bq$.
So we just need to show that $e^\mu_\nu$ acts as zero
on 
$S^{\bq}_{\blambda}$
for all $\blambda \in \MP_d(l)$ 
with $\bi^\bq_\blambda \not\sim \bi$.
Given such a $\blambda$, this is clear from the
definition of the idempotent $e$ in Theorem~\ref{dcp}
unless $\ell(\lambda^{(r)}) \leq \mu_r$ for each $r=1,\dots,l$.
In that case, there is a unique weight $\alpha \in \Col^\mu_d$
mapping to $\blambda$ under the second bijection
from Lemma~\ref{comb}, and the assumption that
$\bi^\bq_\blambda \not\sim \bi$ is equivalent to the statement 
that $\nu(\alpha) \neq \nu$.
Let $N(\alpha)$ be the parabolic Verma module in $\mathcal O^\mu$
of highest weight $(\alpha-\rho)$.
By \cite[Theorem 6.12]{schur}, we have that
$$
\hom_{\mathfrak{g}}(P^\mu \otimes V^{\otimes d}, N(\alpha))
\cong S^\bq_\blambda
$$
as $H_d^\bq$-modules.
Since 
$N(\alpha)$ and $(P^\mu \otimes V^{\otimes d}) e_\nu^\mu$
have different generalized 
central characters,
we have that
$$
e_\nu^\mu 
\hom_{\mathfrak{g}}(P^\mu \otimes V^{\otimes d}, N(\alpha))
=
\hom_{\mathfrak{g}}((P^\mu \otimes V^{\otimes d}) e_\nu^\mu, N(\alpha))
= 0.
$$
Hence $e_\nu^\mu S^\bq_\blambda = 0$ as required.

To deduce the general case,
assume still that $\mu$ is a partition and 
take another composition $\mu' \sim \mu$.
We can find a central element $z_\nu \in Z(\mathfrak{g})$
that acts on 
$P^\mu \otimes V^{\otimes d}$ in the same way
as $e_\nu^\mu$ and on
$P^{\mu'} \otimes V^{\otimes d}$ in the same way
as $e_\nu^{\mu'}$.
In the notation of Lemma~\ref{ms}, 
we have shown that $\rho^\mu(b(\bi))$ coincides with the endomorphism
of $P^\mu \otimes V^{\otimes d}$
defined by left multiplication by
$z_\nu$. We need to prove that $\rho^{\mu'}(b(\bi))$ does too.
This follows because the map $\iota_{\mu,\mu'}$ commutes with the
action of $z_\nu$.
\end{proof}

\begin{Corollary}\label{nearly}
The map $\rho^\mu_\nu:b(\bi) H_d^\bq \rightarrow
\End_{\mathfrak{g}}((P^\mu \otimes V^{\otimes d}) e^\mu_\nu)^{\op}$
induced by the right action of $H_d^\bq$
on $P^\mu \otimes V^{\otimes d}$ is an isomorphism.
Moreover, the center of this algebra is generated by the
endomorphisms $z_1,\dots,z_n$, and is of dimension 
equal to the number of isomorphism classes of irreducible modules
in $\mathcal O^\mu_\nu$.
\end{Corollary}

\begin{proof}
The first statement is immediate from Theorem~\ref{dcp}
since $e^\mu_\nu = b(\bi)$.
The fact that the center of 
$\End_{\mathfrak{g}}((P^\mu \otimes V^{\otimes d}) e^\mu_\nu)^{\op}$
is generated by the endomorphisms $z_1,\dots,z_n$
is immediate from Corollary~\ref{C2}.
The center is isomorphic to $Z(b(\bi) H_d^\bq)$,
which by Theorem~\ref{easiest} is of
dimension equal to the 
size of the set $\{\blambda \in \MP_d(l)\mid \bi^\bq_\blambda = \bi\}$.
By Lemma~\ref{comb}, this is the same as the size of the set
$\Col^\mu_\nu$, that is, the number of isomorphism classes of irreducible
modules in $\mathcal O^\mu_\nu$.
\end{proof}

We need just one more fact, which is a variation on a result of Irving \cite{I}.

\begin{Lemma}\label{sc}
For any $\nu\in Y^\mu_d$, 
the injective hull of any 
module in $\mathcal O^\mu_\nu$ with a parabolic Verma
flag is a finite direct sum of 
direct summands of
$(P^\mu \otimes V^{\otimes d}) e^\mu_\nu$.
\end{Lemma}

\begin{proof}
We claim that every irreducible submodule of a parabolic Verma module
in $\mathcal O^\mu_\nu$ 
embeds into $(P^\mu \otimes V^{\otimes d}) e^\mu_\nu$.
Since $(P^\mu \otimes V^{\otimes d}) e^\mu_\nu$ 
is injective this implies the lemma.
To prove the claim, 
recall that $Y^\mu_d$ is the set that parametrizes the central characters
arising from irreducible constituents of $P^\mu \otimes V^{\otimes d}$.
So the claim follows by \cite[Theorem 4.8]{schur}
in the special case that $\mu$ is actually a partition.
Essentially the same proof as there proves the analogue of this theorem
for arbitrary $\mu$, 
providing one 
replaces the definition of ``standard tableau'' used in \cite{schur}
with the less familiar notion from \cite[(2.2)]{dual}.
\end{proof}

Now finally we consider the commutative diagram
$$
\begin{CD}
&&\:\:\:\,\:\:\:Z(\mathfrak{g}) \\
\\
Z(\mathcal O^\mu_\nu)&@>>f^\mu_\nu>&
Z(\End_{\mathfrak{g}}((P^\mu \otimes V^{\otimes d})e^\mu_\nu)^\op),
\end{CD}
\begin{picture}(0,0)
\put(-176,0){\makebox(0,0){$\swarrow$}}
\put(-176,0){\line(1,1){18}}
\put(-107,0){\makebox(0,0){$\searrow$}}
\put(-107,0){\line(-1,1){18}}
\put(-183,11){\makebox(0,0){$_{m^\mu_\nu}$}}
\put(-102,11){\makebox(0,0){$_{g^\mu_\nu}$}}
\end{picture}
$$
where $m^\mu_\nu, f^\mu_\nu$ and $g^\mu_\nu$ are 
the natural multiplication maps.

\begin{Theorem}\label{there}
In the above diagram, the maps $m^\mu_\nu$ and $g^\mu_\nu$ are 
surjective and
the map $f^\mu_\nu$ is an isomorphism.
Hence, $Z(\mathcal O^\mu_\nu)$ is isomorphic
to $Z(b(\bi) H_d^\bq)$ and is of dimension equal to the number
of isomorphism classes of irreducible module in $\mathcal O^\mu_\nu$.
\end{Theorem}

\begin{proof}
We first prove that $f^\mu_\nu$ is injective.
Suppose we are given a natural transformation
$z \in Z(\mathcal O^\mu_\nu)$ defining
the zero endomorphism of $(P^\mu \otimes V^{\otimes d})e^\mu_\nu$.
To show that $z = 0$, we need to show that $z$ defines the
zero endomorphism of every module $M \in \mathcal O^\mu_\nu$.
Let $P$ be the projective cover of $M$ and $I$ be the injective
hull of $P$. Since $P$ has a parabolic Verma flag by general theory,
Lemma~\ref{sc} implies that $I$ is a finite direct sum of
summands of $(P^\mu \otimes V^{\otimes d})e^\mu_\nu$.
Hence $z$ defines the zero endomorphism of $I$.
Since $P$ embeds into $I$ and surjects onto $M$, we get from this
that $z$ defines the zero endomorphism of $M$ too.
Now to finish the proof of the theorem, 
we know already from Corollary~\ref{nearly}
that $g^\mu_\nu$ is surjective. Hence by the commutativity of the diagram,
$m^\mu_\nu$ and $f^\mu_\nu$ must both be surjetive too.
The remaining statements are immediate from Corollary~\ref{nearly}.
\end{proof}

Finally, we note for any $\nu \in Y^\mu$ that
tensoring with a sufficiently large
power of determinant induces an equivalence between
$\mathcal O^\mu_\nu$ and $\mathcal O^\mu_{\nu'}$ 
for some $\nu' \in Y^\mu_d$ and some $d \geq 0$.
Given this, Theorem 2 from the introduction follows
from Theorem~\ref{there}.


\begin{thebibliography}{AMR}

\bibitem[AS]{AS}
T. Arakawa and T. Suzuki,
Duality between $\mathfrak{sl}_n(\mathbb{C})$ and the degenerate affine Hecke algebra, {\em  J. Algebra} {\bf  209} (1998), 288--304; {\tt q-alg/9710037}.

\bibitem[AMR]{AMR}
S. Ariki, A. Mathas and H. Rui,
Cyclotomic Nazarov-Wenzl algebras, to appear in {\em Nagoya Math. J.};
{\tt math.QA/0506467}.

\bibitem[BG]{BG}
K. Brown and K. Goodearl,
{\em Lectures on algebraic quantum groups},
Birkh\"auser, 2002.

\iffalse
\bibitem[BK1]{BKcrystal}
J. Brundan and A. Kleshchev, 
Hecke-Clifford superalgebras, crystals of type $A_{2\ell}^{(2)}$
and 
modular branching rules for $\widehat{S}_n$, {\em
Represent. Theory} {\bf 5} (2001), 317-403; {\tt math.RT/0103060}. 
\fi

\bibitem[B1]{dual}
J. Brundan,
Dual canonical bases and Kazhdan-Lusztig polynomials,
to appear in {\em J. Algebra};
{\tt math.QA/0509700}.

\bibitem[B2]{Bpre}
J. Brundan,
Symmetric functions, parabolic category $\mathcal O$ and the Springer fiber;
{\tt math.RT/0608235}.

\iffalse
\bibitem[BK1]{rep}
J. Brundan and A. Kleshchev,
Representations of shifted Yangians and finite $W$-algebras,
to appear in {\em Mem. Amer. Math. Soc.};
{\tt math.RT/0508003}.
\fi

\bibitem[BK]{schur}
J. Brundan and A. Kleshchev,
Schur-Weyl duality for higher levels;
{\tt math.RT/0605217}.

\bibitem[D]{D}
V. Drinfeld, 
Degenerate affine Hecke algebras and Yangians, {\em 
Func. Anal. Appl.} {\bf 20} 
(1986), 56--58.

\bibitem[E]{E}
D. Eisenbud,
{\em Commutative algebra with a view towards algebraic geometry},
GTM 150, Springer, 1994.

\bibitem[G]{Gnote}
I. Grojnowski,
Blocks of the cyclotomic Hecke algebra, preprint, 1999.

\bibitem[I]{I}
R. Irving,
Projective modules in the category $\mathcal O_S$: self-duality,
{\em Trans. Amer. Math. Soc.} {\bf 291} (1985), 701--732.

\bibitem[J]{J}
A. Jucys, 
Symmetric polynomials and the center of the symmetric group ring,
{\em Report Math. Phys.} {\bf 5} (1974), 107--112. 

\bibitem[Kh]{Kh}
M. Khovanov, Crossingless matchings and the cohomology of
$(n,n)$ Springer varieties,
{\em Commun. Contemp. Math.} {\bf 6} (2004), 561--577;
{\tt math.QA/0202110}.

\bibitem[K]{Kbook}
A. Kleshchev, {\em
Linear and Projective Representations of Symmetric Groups}, Cam- 
bridge University Press, Cambridge, 2005. 

\bibitem[L]{L}
G. Lusztig,
 Cuspidal local systems and graded Hecke algebras. I. {\em 
Inst. Hautes \'Etudes Sci. Publ. Math.} {\bf 67} (1988), 145--202.

\bibitem[LM]{LM}
S. Lyle and A. Mathas,
The blocks of the affine and cyclotomic Hecke algebras;
{\tt math.RT/0607451}.

\bibitem[Mac]{Mac}
I. Macdonald,
{\em Symmetric functions and Hall polynomials},
Oxford Mathematical Monographs, second edition, OUP, 1995.

\bibitem[MS]{MS}
V. Mazorchuk and C. Stroppel, Projective-injective modules, Serre functors and 
symmetric algebras; {\tt math.RT/0508119}. 

\bibitem[M]{M}
G. Murphy, 
The idempotents of the symmetric group and Nakayama's conjecture, 
{\em J. Algebra} {\bf 81} (1983), 258--265. 

\bibitem[S]{St}
C. Stroppel,
Perverse sheaves on Grassmannians, Springer fibres and 
Khovanov homology; {\tt math.RT/0608234}.

\bibitem[W]{wang}
W. Wang,
Vertex algebras and the class algebras of wreath products,
{\em Proc. London Math. Soc.} {\bf 88} (2004), 381--404; {\tt math.QA/0203004}.
\end{thebibliography}
\end{document}